\documentclass[leqno,12pt]{article}

\usepackage{amssymb}
\usepackage{euscript}
\usepackage[dvips]{graphicx}
\usepackage{flafter}
\usepackage{pstricks}

\usepackage{verbatim}

\usepackage{amsmath}

\usepackage{color}

\usepackage{mathrsfs} 

\usepackage{amsthm}

\evensidemargin\oddsidemargin

\begin{document}

\baselineskip=18pt
\setcounter{page}{1}

\renewcommand{\theequation}{\thesection.\arabic{equation}}
\newtheorem{theorem}{Theorem}[section]
\newtheorem{lemma}[theorem]{Lemma}
\newtheorem{definition}[theorem]{Definition}
\newtheorem{proposition}[theorem]{Proposition}
\newtheorem{corollary}[theorem]{Corollary}
\newtheorem{fact}[theorem]{Fact}
\newtheorem{problem}[theorem]{Problem}
\newtheorem{conjecture}[theorem]{Conjecture}
\newtheorem{claim}[theorem]{Claim}

\theoremstyle{definition} 
\newtheorem{remark}[theorem]{Remark}

\newcommand{\eqnsection}{
\renewcommand{\theequation}{\thesection.\arabic{equation}}
    \makeatletter
    \csname  @addtoreset\endcsname{equation}{section}
    \makeatother}
\eqnsection


\def\r{{\mathbb R}}
\def\e{{\mathbb E}}
\def\p{{\mathbb P}}
\def\P{{\bf P}}
\def\E{{\bf E}}
\def\Q{{\bf Q}}
\def\z{{\mathbb Z}}
\def\N{{\mathbb N}}
\def\T{{\mathbb T}}
\def\G{G}

\def\ee{\mathrm{e}}
\def\d{\, \mathrm{d}}
\def\dtv{\mathrm{d_{tv}}}



\vglue50pt

\centerline{\Large\bf The slow regime of randomly biased walks on trees}

{
\let\thefootnote\relax\footnotetext{\scriptsize Partly supported by ANR project MEMEMO2 (2010-BLAN-0125).}
}

\bigskip
\bigskip

\centerline{by}

\medskip

\centerline{Yueyun Hu\footnote{\scriptsize LAGA, Universit\'e Paris XIII, 99 avenue J-B Cl\'ement, F-93430 Villetaneuse, France, {\tt yueyun@math.univ-paris13.fr}} and Zhan Shi\footnote{\scriptsize LPMA, Universit\'e Paris VI, 4 place Jussieu, F-75252 Paris Cedex 05, France, {\tt zhan.shi@upmc.fr}}}

\medskip

\centerline{\it Universit\'e Paris XIII  \& Universit\'e Paris VI}



\bigskip
\bigskip
\bigskip

{\leftskip=2truecm \rightskip=2truecm \baselineskip=15pt \small

\noindent{\slshape\bfseries Summary.} We are interested in the randomly biased random walk on the supercritical Galton--Watson tree. Our attention is focused on a slow regime when the biased random walk $(X_n)$ is null recurrent, making a maximal displacement of order of magnitude $(\log n)^3$ in the first $n$ steps.   We study the localization problem of $X_n$ and prove that the quenched law of $X_n$ can be approximated by a certain invariant probability   depending on $n$ and the random environment. As a consequence, we establish that    upon the survival  of the system, $\frac{|X_n|}{(\log n)^2}$ converges in law to some non-degenerate limit on $(0, \infty)$ whose law is explicitly computed. 
\bigskip


\bigskip


} 

\bigskip
\bigskip

\section{Introduction}
   \label{s:intro}

$\phantom{aob}$Let $\T$ be a supercritical Galton--Watson tree rooted at $\varnothing$, so it survives with positive probability. For any pair of vertices $x$ and $y$ of $\T$, we say $x\sim y$ if $x$ is either a child, or the parent, of $y$. Let $\omega := (\omega(x), \, x\in \T)$ be a sequence of vectors; for each vertex $x\in \T$, $\omega(x) := (\omega(x, \, y), \, y\in \T)$ is such that $\omega(x, \, y) \ge 0$ for all $y\in \T$ and that $\sum_{y\in \T} \omega(x, \, y) =1$. We assume that for each pair of vertices $x$ and $y$, $\omega(x, \, y)>0$ if and only if $y \sim x$.

For given $\omega$, let $(X_n, \, n\ge 0)$ be a random walk on $\T$ with transition probabilities $\omega$, i.e., a $\T$-valued Markov chain, started at $X_0 = \varnothing$, such that
$$
P_\omega \{ X_{n+1} = y \, | \, X_n =x \} 
= 
\omega(x, \, y).
$$

For any vertex $x\in \T \backslash \{ \varnothing\}$, let ${\buildrel \leftarrow \over x}$ be its parent, and let $(x^{(1)}, \cdots, x^{(N(x))})$ be its children, where $N(x) \ge 0$ is the number of children of $x$. Define $A(x) := (A_i(x), \, 1\le i\le N(x))$ by
\begin{equation}
    A_i(x) 
    := 
    \frac{\omega(x, \, x^{(i)})}{\omega(x, \, {\buildrel\leftarrow \over x})},
    \qquad 1\le i\le N(x) \, .
    \label{A}
\end{equation}

\noindent A special example is when $A_i(x) = \lambda$ for all $x\in \T\backslash \{ \varnothing \}$ and all $1\le i\le N(x)$, where $\lambda$ is a finite and positive constant: the random walk $(X_n)$ is then the $\lambda$-biased random walk on $\T$ introduced and studied in depth by Lyons~\cite{lyons90}--\cite{lyons92}, Lyons, Pemantle and Peres~\cite{lyons-pemantle-peres-ergodic}--\cite{lyons-pemantle-peres96}.
In particular, if $A_i(x) = 1$, $\forall x$, $\forall i$, we get the simple random walk on $\T$. 

It is known that when the transition probabilities are {\it random} --- the resulting random walk $(X_n)$ is then a random walk in random environment ---, the walk possesses a regime of {\it slow movement}. We are interested in this slow movement in this paper. 

In the language of Neveu~\cite{neveu86}, $(\T, \, \omega)$ is a marked tree. Note that $A(x)$, $x\in \T \backslash \{ \varnothing\}$ depends entirely on the marked tree. We assume, from now on, that $A(x)$, $x\in \T \backslash \{ \varnothing\}$, are i.i.d., and write $A=(A_1, \cdots, A_N)$ for a generic random vector having the law of $A(x)$ (for any $x\in \T \backslash \{ \varnothing\}$). We mention that the dimension $N \ge 0$ of $A$ is random, and is governed by the law of reproduction of $\T$. We use $\P$ to denote the probability with respect to the environment, and $\p := \P \otimes P_\omega$ the annealed probability, i.e., $\p(\, \cdot \, ) := \int P_\omega( \, \cdot \, ) \, \P(\! \d \omega)$. 

Throughout the paper, we assume
\begin{equation}
    \E \Big( \sum_{i=1}^N A_i \Big) 
    = 
    1,
    \qquad
    \E \Big( \sum_{i=1}^N A_i \log A_i \Big) 
    =
    0 \, .
    \label{cond-hab-A}
\end{equation}

\noindent In the language of branching random walks (see Section \ref{s:lt}), (\ref{cond-hab-A}) refers to the ``boundary case"; in this case, the biased walks produce some unusual phenomena that have still been beyond good understanding. We also assume the following integrability condition: there exists $\delta>0$ such that
\begin{equation}
    \E \Big( \sum_{i=1}^N A_i^{1+\delta} \Big)
    +
    \E \Big( \sum_{i=1}^N A_i^{-\delta} \Big) 
    +
    \E (N^{1+\delta}) 
    <
    \infty .
    \label{integrability-assumption-A}
\end{equation}

\noindent Lyons and Pemantle~\cite{lyons-pemantle} established a recurrence vs.\ transience criterion for random walks on general trees; applied to the special setting of Galton--Watson trees, it says that (\ref{cond-hab-A}) ensures that the biased walk $(X_n)$ is $\p$-a.s.\ recurrent. Menshikov and Petritis~\cite{menshikov-petritis} gave another proof of the recurrence by means of Mandelbrot's multiplicative cascades, assuming some additional integrability condition. The proofs of the recurrence in both \cite{lyons-pemantle} and \cite{menshikov-petritis} required an extra exchangeability assumption on $(A_1, \, \cdots, \, A_N)$, which turned out to be superfluous, and was removed by Faraud~\cite{faraud}, who furthermore proved that $(X_n)$ is null recurrent under (\ref{cond-hab-A}).

Introduced by Lyons and Pemantle~\cite{lyons-pemantle} as an extension of deterministically biased walks studied in Lyons~\cite{lyons90}-\cite{lyons92}, randomly biased walks on trees have received much research interest. Deep results were obtained by Lyons, Pemantle and Peres~\cite{lyons-pemantle-peres-ergodic} and \cite{lyons-pemantle-peres96}, who also raised further open problems. Often motivated by these results and problems, both the transient regimes (\cite{elie}, \cite{elie2}) 
and the recurrent regimes (\cite{andreoletti-debs1}, \cite{andreoletti-debs2}, \cite{faraud}, \cite{gyzbiased}, \cite{yztree},  \cite{yzslow}) have been under intensive study for these walks. For a general account of biased walks on trees, we refer to \cite{lyons-peres}, \cite{peres} and \cite{stf}.

We add a special vertex, denoted by ${\buildrel \leftarrow \over \varnothing}$, which is the parent of $\varnothing$, and assume that $(\omega(\varnothing, \, y), \, |y|=1 \; \mathrm{or} \; y={\buildrel \leftarrow \over \varnothing})$ is independent of other random vectors $(\omega(x, \, y), \, y\sim x)$ for $x\in \T \backslash \{ \varnothing\}$, having the same distribution as any of these random vectors; whenever the biased walk $(X_i)$ hits ${\buildrel \leftarrow \over \varnothing}$, it comes back to $\varnothing$ in the next step. [However, ${\buildrel \leftarrow \over \varnothing}$ is not considered as a vertex of $\T$; so, for example, $\sum_{x\in \T} f(x)$ does not contain the term $f({\buildrel \leftarrow \over \varnothing})$.] This makes the presentation of our model more pleasant, since the family of i.i.d.\ random vectors $A(x)$ also includes the element $A(\varnothing)$ from now on.

It was proved in \cite{gyzbiased} that under \eqref{cond-hab-A} and \eqref{integrability-assumption-A}, almost surely upon the survival of the system, 
$$
\lim_{n\to\infty} \frac{1}{(\log n)^3} \max_{0\le i\le n} |X_i| 
=
\frac{8}{3 \pi^2 \sigma^2} \, ,
$$ 

\noindent where
\begin{equation}
    \label{sigma2}   
    \sigma^2
    := 
    \E \Big( \sum_{i=1}^N \, A_i  (\log A_i)^2  \Big) \in (0, \,\infty). 
\end{equation}

We are interested in the typical size of $|X_n|$; a natural question is to find a deterministic sequence $a_n\to \infty$ such that $\frac{|X_n|}{a_n}$ converges in law to some non-degenerate limit. In dimension 1 (which would be an informal analogue of the case $N(x)=1$ for all $x$), the slow movement was discovered by Sinai~\cite{sinai} who showed that $\frac{X_n}{(\log n)^2}$ converges weakly to a non-degenerate limit under the annealed measure. More precisely, Sinai~\cite{sinai} developed the seminal ``method of valley" to localize the walk around the bottom of a certain Brownian valley with high probability. This method, however, seems hopeless to be directly adapted to the biased walk on trees. Observe that in terms of the invariant measure, we can interpret Sinai's method of valley as the approximation of the law of the walk by a certain invariant measure whose mass is concentrated at the neighbourhood of the bottom. Our main result, stated as Theorem \ref{t:cv} below, asserts that upon the survival of the system, the (quenched) finitely-dimensional distribution of the biased walk can be approximated by the product measure of some invariant probability measures. A consequence of this result is that under \eqref{cond-hab-A} and \eqref{integrability-assumption-A}, for all $x>0$,
$$ 
 \lim_{n\to \infty} \p \Big( \frac{\sigma^2\, |X_n|}{(\log n)^2} \le x \, \Big| \, \mathrm{survival} \Big) 
 =
 \int_0^x \frac{1}{(2\pi y)^{1/2}} \, 
 \P \Big( \eta \le \frac{1}{y^{1/2}} \Big) \d y \, , 
$$ 

\noindent where $\sigma$ is the constant in \eqref{sigma2}, and $\eta := \sup_{s\in [0, \, 1]} [\, \overline{\mathfrak{m}}(s) - \mathfrak{m}(s)]$. Here, $(\mathfrak{m}(s), \, s\in [0, \, 1])$ is a standard Brownian meander,\footnote{Recall that the standard Brownian meander can be realized as follows: $\mathfrak{m}(s) := \frac{|B(\mathfrak{g}+s(1-\mathfrak{g}))|}{(1-\mathfrak{g})^{1/2}}$, $s\in [0, \, 1]$, where $(B(t), \, t\in [0, \, 1])$ is a standard Brownian motion, with $\mathfrak{g} := \sup\{ t\le 1: \, B(t) =0\}$.} and $\overline{\mathfrak{m}}(s) := \sup_{u\in [0, \, s]} \mathfrak{m}(u)$.

We mention that $\int_0^\infty \frac{1}{(2\pi y)^{1/2}} \, \P ( \eta \le \frac{1}{y^{1/2}} ) \d y  =1$ because $\E( \frac{1}{\eta})= (\frac{\pi}{2})^{1/2}$, see \cite{Yor}.

In the next section, we give a precise statement of Theorem \ref{t:cv}, as well as an outline of its proof.


\section{Random potential, and statement of results}
\label{s:lt}

$\phantom{aob}$The movement of the biased random walk $(X_n)$ depends strongly on the random environment $\omega$. It turns out to be more convenient to quantify the influence of the random environment via the random {\bf potential}, which we define by
$V(\varnothing):=0$ and
\begin{equation}
    V(x) 
    := 
    -
    \sum_{y\in \, ]\!] \varnothing,\, x]\!]}
    \log \,
    \frac{\omega({\buildrel \leftarrow \over y},
    \, y)}{\omega({\buildrel \leftarrow \over y}, \,
    {\buildrel \Leftarrow \over y})},
    \qquad x\in \T\backslash\{ \varnothing\} \, ,
    \label{V}
\end{equation}

\noindent where ${\buildrel \Leftarrow \over y}$ is the parent of ${\buildrel \leftarrow \over y}$, and $\, ]\!] \varnothing, \, x]\!] := [\![ \varnothing, \, x]\!] \backslash \{ \varnothing\}$, with $[\![ \varnothing, \, x]\!]$ denoting the set of vertices (including $x$ and $\varnothing$) on the unique shortest path connecting $\varnothing$ to $x$. Throughout the paper, we use $x_i$ (for $0\le i\le |x|$) to denote the ancestor of $x$ in the $i$-th generation; in particular, $x_0 = \varnothing$ and $x_{|x|} =x$. As such, the potential $V$ in (\ref{V}) can also be written as
$$
V(x) 
= 
- \sum_{i=0}^{|x|-1} \log \,
\frac{\omega(x_i, \, x_{i+1})}{\omega(x_i, \, x_{i-1})},
\qquad x\in \T\backslash\{ \varnothing\} \, .
\qquad
(x_{-1} := {\buildrel \leftarrow \over \varnothing})
$$

\noindent The random potential process $(V(x), \, x\in \T)$ is a branching random walk, in the usual sense of Biggins~\cite{biggins77}. There exists an obvious bijection between the random environment $\omega$ and the random potential $V$.

In terms of the random potential, assumptions (\ref{cond-hab-A}) and (\ref{integrability-assumption-A}) become, respectively,
\begin{equation}
    \E \Big( \sum_{x: \, |x|=1} \ee^{-V(x)} \Big) =1,
    \qquad
    \E \Big( \sum_{x: \, |x|=1} V(x)\, \ee^{-V(x)} \Big) =0 ,
    \label{cond-hab}
\end{equation}

\noindent and 
\begin{equation}
    \E \Big( \sum_{x: \, |x|=1} \ee^{-(1+\delta)V(x)} \Big)
    +
    \E \Big( \sum_{x: \, |x|=1} \ee^{ \delta V(x)} \Big) 
    +
    \E \Big[ \Big( \sum_{x: \, |x|=1} 1\Big)^{1+\delta} \, \Big] <\infty \, .
    \label{integrability-assumption}
\end{equation}

\noindent We refer from now on to (\ref{cond-hab}) or (\ref{integrability-assumption}) instead of to (\ref{cond-hab-A}) or (\ref{integrability-assumption-A}). In the language of branching random walks, (\ref{cond-hab}) corresponds to the ``boundary case" (Biggins and Kyprianou~\cite{biggins-kyprianou05}). The branching random walk in this case is known, under some additional integrability assumptions, to have some highly non-trivial universality properties. 

We are often interested in properties upon the system's non-extinction, so let us introduce
\begin{eqnarray*}
    \P^*(\, \cdot \, )
 &:=&\P( \, \cdot \, | \, \hbox{non-extinction}) ,
    \\
    \p^*(\, \cdot \, )
 &:=&\p( \, \cdot \, | \, \hbox{non-extinction}) \, .
\end{eqnarray*}


Let us define a symmetrized version of the potential:
\begin{equation}
    U(x) 
    := 
    V(x) - \log (\frac{1}{\omega(x, \, {\buildrel \leftarrow \over x})}) \, ,
    \qquad
    x\in \T \, .
    \label{U}
\end{equation}

\noindent We call $U$ the {\bf symmetrized potential}, and use frequently the following relation between $U$ and $V$:
\begin{equation}
    \ee^{-U(x)}
    = 
    \frac{1}{\omega(x, \, {\buildrel \leftarrow \over x})} \, \ee^{-V(x)}  
    =
    \ee^{-V(x)} + \sum_{y\in \T: \, {\buildrel \leftarrow \over y} =x} \ee^{-V(y)} ,
    \qquad
    x\in \T \, .
    \label{Theta}
\end{equation}

We now introduce a pair of fundamental martingales associated with the potential $V$.
Assumption (\ref{cond-hab}) immediately implies that $(W_n, \, n\ge 0)$ and $(D_n, \, n\ge 0)$ are martingales under $\P$, where
\begin{eqnarray}
    W_n 
 &:=& \sum_{x: \, |x|=n} \ee^{-V(x)}, 
    \label{Wn}
    \\
    D_n 
 &:=& \sum_{x: \, |x|=n} V(x) \ee^{-V(x)}, 
    \qquad
    n\ge 0 \, ,
    \label{Dn}
\end{eqnarray}

\noindent In the literature, $(W_n)$ is referred to as an {\it additive martingale}, and $(D_n)$ a {\it derivative martingale}. Since $(W_n)$ is a non-negative martingale, it converges $\P$-a.s.\ to a finite limit; under assumption (\ref{cond-hab}), this limit is known (Biggins~\cite{biggins-mart-cvg}, Lyons~\cite{lyons}) to be $0$:
\begin{equation}
    W_n
    \to
    0,
    \qquad
    \hbox{\rm $\P^*$-a.s.}
    \label{W->0}
\end{equation}

\noindent [We will see in (\ref{ezratio}) the rate of convergence.] In view of (\ref{Theta}), this yields 
\begin{equation}
    \inf_{x: \, |x|=n} U(x) 
    \to 
    \infty ,
    \qquad
    \hbox{\rm $\P^*$-a.s.}
    \label{U->infty}
\end{equation}

\noindent For the derivative martingale $(D_n)$, it is known (Biggins and Kyprianou~\cite{biggins-kyprianou04}, A\"id\'ekon~\cite{elie-min}) that (\ref{integrability-assumption}) is ``slightly more than" sufficient to ensure that $D_n$ converges $\P$-a.s.\ to a limit, denoted by $D_\infty$, and that 
$$
D_\infty >0,
\qquad
\hbox{\rm $\P^*$-a.s.}
$$

\noindent For an optimal condition (of $L \log  L$-type) for the positivity of $D_\infty$, see the recent work of Chen~\cite{chen}. The two martingales $(D_n)$ and $(W_n)$ are asymptotically related; see Section \ref{s:proof-Zn}.

The basic idea  is to add a reflecting barrier at (notation: $]\! ]\varnothing, \, x[\![ \; := \; ]\! ]\varnothing, \, x]\!] \backslash \{ x\}$)
\begin{equation}
    \mathscr{L}_r
    := 
    \Big\{ x: \, \sum_{z\in \, ]\! ]  \varnothing, \, x]\! ]} \ee^{V(z)-V(x)} > r , \;
    \sum_{z\in \, ]\! ] \varnothing, \, y]\! ]} \ee^{V(z)-V(y)} \le r, \; 
    \forall y\in \, ]\! ]\varnothing, \, x[\![ \Big\} ,
    \label{gamma}
\end{equation}

\noindent where $r>1$ is a parameter.\footnote{That is, each time the biased walk $(X_i)$ hits any vertex $x\in \mathscr{L}_r$, it moves back to ${\buildrel \leftarrow \over x}$ in the next step.} We mention that $\mathscr{L}_r$ does not necessarily separate $\varnothing$ from infinity: our assumptions (\ref{cond-hab}) and (\ref{integrability-assumption}) do not exclude the existence of $r>1$ and a sequence of vertices $x_0:= \varnothing < x_1 <x_2 \cdots$ with $|x_i| =i$, $i\ge 0$, such that $\sum_{i=1}^n \ee^{V(x_i)-V(x_n)} \le r$ for all $n\ge 1$.

If $r=r(n) := \frac{n}{(\log n)^\gamma}$ with $\gamma<1$, then we will see from Lemma \ref{l:ligne-arret} that with $\p^*$-probability going to 1 (for $n\to \infty$), the biased walk does not hit any vertex in $\mathscr{L}_r$ in the first $n$ steps.\footnote{Actually $\gamma<2$ will do the job (by Theorem \ref{t:ligne-arret}). However, in Section \ref{s:proof-p:x=0}, when we start proving our main results, only Lemma \ref{l:ligne-arret} is available, which says that $\gamma<1$ suffices. The proof of Theorem \ref{t:ligne-arret} comes afterwards, in Section \ref{s:proof-t:ligne-arret}.} As such, it makes no significant difference if we add a reflecting barrier at $\mathscr{L}_r$. An advantage, with the presence of the reflecting barrier at $\mathscr{L}_r$, for any $r>1$, is that the biased walk becomes {\it positive recurrent} under the quenched probability $P_\omega$, and its invariant probability $\pi_r$ is as follows: $\pi_r({\buildrel \leftarrow \over \varnothing}) := \frac{1}{Z_r}$, and for $x\in \T$,
\begin{equation}
    \pi_r(x)
    :=
    \begin{cases}
    \frac{1}{Z_r} \, \ee^{-U(x)} \, ,
\qquad \hbox{\rm if $x<\mathscr{L}_r$} \, ,
    \\
    \\
    \frac{1}{Z_r} \, \ee^{-V(x)} \, , \qquad \hbox{\rm if $x\in \mathscr{L}_r$} \, ,
    \end{cases}
    \label{pi_n}
\end{equation}

\noindent where $Z_r$ is the normalizing factor:\footnote{By $x<\mathscr{L}_r$, we mean $\sum_{z\in \, ]\! ] \varnothing, \, y]\! ]} \ee^{V(z)-V(y)} \le r$ for all vertex $y\in \, ]\! ]\varnothing, \, x]\!]$.}
\begin{equation}
    Z_r
    :=
    1+ \sum_{x\in \T: \, x<\mathscr{L}_r} \ee^{-U(x)} 
    + \sum_{x\in \mathscr{L}_r} \ee^{-V(x)} .
    \label{Zn} 
\end{equation}
    
\noindent We extend the definition of $\pi_r$ to the whole tree $\T$ by letting  $\pi_r(x):= 0$ if  neither $x<\mathscr{L}_r$ nor $x \in \mathscr{L}_r$. 

Due to the periodicity of the walk $(X_i)$, we divide the tree $\T$ into $\T^{(\mathrm{even})}$ and $\T^{(\mathrm{odd})}$ with $$ \T^{(\mathrm{even})}:= \{ x \in \T: |x| \hbox{ is even}\}, \qquad \T^{(\mathrm{odd})}:= \{ x \in \T: |x| \hbox{ is odd}\}.$$ 

Depending on the parity of $n$,  the law of $X_n$ (starting from $\varnothing$) is supported either by  $ \T^{(\mathrm{even})}$ or by $\T^{(\mathrm{odd})} \cup\{{\buildrel \leftarrow \over \varnothing}\}$. Note that $\pi_r ( \T^{(\mathrm{even})}) = \pi_r ( \T^{(\mathrm{odd})} \cup\{{\buildrel \leftarrow \over \varnothing}\})= \frac12$ as $\pi_r (\cdot)$ is the invariant probability measure of a finite Markov chain of period $2$. We define a new probability measure: for any $r>1$,
\begin{equation}
    \widetilde \pi_r(\, \cdot\, )
    :=     
    \begin{cases}
    2 \pi_r (\cdot) \, {\bf 1}_{ \T^{(\mathrm{even})}}(\cdot) \, ,
    & \hbox{\rm if $\lfloor r\rfloor$ is even} ,
    \\
    \\
   2  \pi_r (\, \cdot\, ) \, {\bf 1}_{ \T^{(\mathrm{odd})} \cup\{{\buildrel \leftarrow \over \varnothing}\}} (\cdot) \, ,
    & \hbox{\rm if $\lfloor r\rfloor$ is odd} .
    \end{cases}
    \label{pi_nbis}
\end{equation}
 
For any pair of probability measures $\mu$ and $\nu$ on $\T\cup\{{\buildrel \leftarrow \over \varnothing}\}$, we denote by $\dtv(\mu, \nu)$ the distance in total variation: 
$$
\dtv(\mu, \nu)
:= 
\frac12 \sum_{x \in \T\cup\{{\buildrel \leftarrow \over \varnothing}\}} |\mu(x) - \nu(x)| \, .
$$  

The main result of the paper is as follows.

\medskip

\begin{theorem}
 \label{t:cv}
 
 Assume $(\ref{cond-hab})$ and $(\ref{integrability-assumption})$. Then
 $$ 
 \dtv\Big( P_\omega\{ X_n \in \bullet\} \, , \; \widetilde \pi_n \Big) 
 \to 
 0, 
 \qquad\hbox{in $\P^*$-probability} \, .
 $$
 More generally, for any $\kappa \ge 1$ and $0<t_1<t_2<\cdots<t_\kappa \le 1$,  
 $$ 
 \dtv\Big( P_\omega\{ (X_{\lfloor t_1 n\rfloor}, \, \cdots \, , X_{\lfloor t_\kappa n\rfloor}) \in \bullet\} \, , \; \bigotimes_{i=1}^\kappa \widetilde \pi_{t_i n} \Big) 
 \to 
 0, 
 \qquad\hbox{in $\P^*$-probability} \, .
 $$
\end{theorem}

\medskip

As such, $X_{\lfloor t_i n\rfloor}$, $1\le i \le \kappa$, are asymptotically independent under $P_\omega$. In particular, no aging phenomenon is possible in the scale of linear time. 

Let us mention that in Theorem \ref{t:cv}, the dependence of $\widetilde \pi_{t_i n}$ on $t_i$ is rather weak. As Lemma \ref{l:pinpim} below shows, $\dtv ( \pi_{t_i n}, \, \pi_n) \to 0$ in $\P^*$-probability, so asymptotically, the influence of $t_i$ on $\widetilde \pi_{t_i n}$ shows up only via the parity of $\lfloor t_i n\rfloor$. 

\medskip

\begin{lemma} 
 \label{l:pinpim}    

 For any $a\ge 0$, as $r\to \infty$, 
 $$ 
 \sup_{u\in [\frac{r}{(\log r)^a}, \, r]} \dtv (\pi_r, \, \pi_u) 
 \to 
 0, \qquad  \hbox{in $\P^*$-probability} . 
 $$ 

\end{lemma}

\medskip

Theorem \ref{t:cv} has the following interesting consequence concerning distance between $X_n$ and $\varnothing$.

\medskip

\begin{corollary}
 \label{c:cv}  
 
 Assume $(\ref{cond-hab})$ and $(\ref{integrability-assumption})$. Fix  $\kappa\ge1$ and $0<t_1<t_2<\cdots<t_\kappa \le1$. Under $\p^*$, $\frac{\sigma^2}{(\log n)^2}\, |X_{\lfloor t_i n\rfloor}|$, $1\le i \le \kappa$, are asymptotically independent and converge in law to a common non-degenerate limit on $(0, \, \infty)$ whose density is given by  
 $$
 \frac{1}{(2\pi x)^{1/2}} \, \P \Big( \eta \le \frac{1}{x^{1/2}} \Big) \, {\bf 1}_{\{x >0\}} \, , 
 $$ 
 where $\sigma^2\in (0, \, \infty)$ is the constant in $(\ref{sigma2})$, and $\eta := \sup_{s\in [0, \, 1]} [\, \overline{\mathfrak{m}}(s) - \mathfrak{m}(s)]$. Here, $(\mathfrak{m}(s), \, s\in [0, \, 1])$ is a standard Brownian meander, and $\overline{\mathfrak{m}}(s) := \sup_{u\in [0, \, s]} \mathfrak{m}(u)$.

\end{corollary}

\medskip

The distribution of $\eta$ is easily seen to be absolutely continuous (Section \ref{s:proof-Zn}), and can be characterised using a result of Lehoczky~\cite{lehoczky}. For more discussions, see \cite{Yor}. Very recently, Pitman~\cite{pitman} has succeeded in determining the law of $\eta$ using a relation between the Brownian meander and the Brownian bridge established by Biane and Yor~\cite{biane-yor}: $\eta$ has the Kolmogorov--Smirnov distribution:
$$
\P(\eta\le x)
=
\sum_{k = -\infty}^\infty (-1)^k \ee^{-2 k^2 x^2}
=
\frac{(2\pi)^{1/2}}{x} \sum_{j=0}^\infty \exp\Big( - \frac{(2j+1)^2 \pi^2}{8x^2} \Big) ,
\qquad
x>0.
$$

Theorem \ref{t:cv} is proved by means of two intermediate estimates, stated below as Propositions \ref{p:x=0} and \ref{p:local}. The first proposition estimates the local time at the root $\varnothing$, whereas the second concerns the local limit probability of the biased walk.  

For any vertex $x\in \T$, let us define
\begin{equation}
    L_n(x)
    :=
    \sum_{i=1}^n \, {\bf 1}_{\{ X_i =x\} } \, ,
    \qquad
    n\ge 1 \, ,
    \label{local_time}
\end{equation}

\noindent which is the (site) local time of the biased walk at $x$.  

\medskip

\begin{proposition}
\label{p:x=0}

 Assume $(\ref{cond-hab})$ and $(\ref{integrability-assumption})$. 
 For any $\varepsilon>0$,
\begin{equation}\label{localtimeproba}
 P_\omega 
 \Big\{ \, \Big| \frac{L_n(\varnothing)}{\frac{n}{\log n}} - \frac{\sigma^2}{4 D_\infty} \, \ee^{-U(\varnothing)}\Big| > \varepsilon \Big\}
 \to 
 0, 
 \qquad\hbox{in $\P^*$-probability} \, .
\end{equation}   Moreover,  \begin{equation}\label{localtimemoment}  E_\omega 
 \Big( \frac{L_n(\varnothing)}{\frac{n}{\log n}} \Big) \to  \frac{\sigma^2}{4 D_\infty} \, \ee^{-U(\varnothing)},  \qquad\hbox{in $\P^*$-probability} \, ,
\end{equation}
\end{proposition}

\medskip

 \begin{proposition}   
\label{p:local}
 Assume $(\ref{cond-hab})$ and $(\ref{integrability-assumption})$.   As $n \to \infty$ along even numbers,
 $$ 
 ( \log n) P_\omega(X_n=\varnothing)  \to \frac{\sigma^2}{2 D_\infty} \, \ee^{-U(\varnothing)},  \qquad\hbox{in $\P^*$-probability} \,.
 $$ 
\end{proposition}
 
\medskip

We now say a few words about the proof. It turns out that the partition function $Z_r$ has a simpler expression. Let $\mathscr{L}_r$ be as in \eqref{gamma}. Define
\begin{equation}
    Y_r
    :=
    \sum_{x\in \T: \, x \le \mathscr{L}_r} \ee^{-V(x)} \, ,
    \label{Yn}
\end{equation}

\noindent with the obvious notation $x\le \mathscr{L}_r$ meaning $x<\mathscr{L}_r$ or $x \in \mathscr{L}_r$. 

\medskip

\begin{lemma}
\label{l:Zn=2Yn}
 
  Let $Y_r$ and $Z_r$ be as in $(\ref{Yn})$ and $(\ref{Zn})$, respectively. Then $Z_r = 2Y_r$, for all $r>1$. 

\end{lemma}

\medskip

\noindent {\it Proof.} If $x\in \T$ is such that $x<\mathscr{L}_r$, we have $Z_r \, \pi_r(x) = \ee^{-U(x)} = \ee^{-V(x)} + \sum_{y\in \T: \, {\buildrel \leftarrow \over y} =x} \ee^{-V(y)}$. Therefore,
\begin{eqnarray*}
    \sum_{x<\mathscr{L}_r} Z_r \, \pi_r(x)
 &=&\sum_{x<\mathscr{L}_r}\ee^{-V(x)}
    +
    \sum_{x<\mathscr{L}_r} \sum_{y\in \T: \, {\buildrel \leftarrow \over y} =x} \ee^{-V(y)}
    \\
 &=&\sum_{x<\mathscr{L}_r}\ee^{-V(x)}
    +
    \sum_{y\in \T: \, \varnothing < y \le \mathscr{L}_r} \ee^{-V(y)} \, ,
\end{eqnarray*}

\noindent which is $\sum_{x<\mathscr{L}_r}\ee^{-V(x)} + \sum_{y\le\mathscr{L}_r}\ee^{-V(y)} -\ee^{-V(\varnothing)}$. Hence
$$
    \sum_{x<\mathscr{L}_r} Z_r \, \pi_r(x)
    =
    2\sum_{x\le\mathscr{L}_r}\ee^{-V(x)} - \sum_{x\in\mathscr{L}_r}\ee^{-V(x)} -1 \, .
$$

\noindent Since $\pi$ is a probability measure, we have $\pi_r({\buildrel \leftarrow \over \varnothing}) + \sum_{x<\mathscr{L}_r} \pi_r(x) + \sum_{x \in \mathscr{L}_r} \pi_r(x) =1$, so
\begin{eqnarray*}
    Z_r
 &=& Z_r \, \pi_r({\buildrel \leftarrow \over \varnothing}) + \sum_{x<\mathscr{L}_r} Z_r \, \pi_r(x) + \sum_{x \in \mathscr{L}_r} Z_r \, \pi_r(x)
    \\
 &=&1
    +
    \Big[2\sum_{x\le\mathscr{L}_r}\ee^{-V(x)} - \sum_{x\in\mathscr{L}_r}\ee^{-V(x)} -1\Big]
    +
    \sum_{x \in \mathscr{L}_r} \ee^{-V(x)} \, ,
\end{eqnarray*}

\noindent which is $2\sum_{x\le \mathscr{L}_r}\ee^{-V(x)}$. Lemma \ref{l:Zn=2Yn} is proved.\hfill$\Box$

\bigskip

So $Y_r$ is half the partition function under the invariant measure. The following theorem, which plays an important role in the proof of Proposition \ref{p:x=0} and Theorem \ref{t:cv}, describes the asymptotics of $Y_r$.

\medskip

\begin{theorem}
\label{t:Yn}

 Assume $(\ref{cond-hab})$ and $(\ref{integrability-assumption})$. Let $Y_r$ be as in $(\ref{Yn})$. We have
 $$
    \lim_{r\to \infty} \, \frac{Y_r}{\log r}
    =
    \frac{2}{\sigma^2} \, D_\infty \, ,
    \qquad\hbox{in $\P^*$-probability} \, ,
 $$
 where $\sigma^2\in (0, \, \infty)$ is the constant in $(\ref{sigma2})$, and $D_\infty$ the $\P^*$-almost sure positive limit of the derivative martingale $(D_n)$ in $(\ref{Dn})$. As a consequence,
 $$
    \lim_{r\to \infty} \, \frac{Z_r}{\log r}
    =
    \frac{4}{\sigma^2} \, D_\infty \, ,
    \qquad
    \lim_{r\to \infty} \, (\log r) \, \pi_r(\varnothing)
    =
    \frac{\sigma^2}{4D_\infty} \, \ee^{-U(\varnothing)} \, ,    
    \qquad\hbox{in $\P^*$-probability} \, .
 $$ 

\end{theorem}

\medskip

Finally, the following general estimate allows us to justify the presence of a barrier at $\mathscr{L}_r$.

\medskip

\begin{theorem}
\label{t:ligne-arret}

 Assume $(\ref{cond-hab})$ and $(\ref{integrability-assumption})$. Let $(a_n)$ be a deterministic sequence of positive real numbers such that $\,\lim_{n\to \infty} \frac{a_n}{(\log n)^2} = 0$, then
 $$
 \lim_{n\to \infty} \, 
 \p \Big( \bigcup_{i=1}^n\{ X_i \in \mathscr{L}_{r_n} \} \Big)
 =
 0 \, ,
 $$
 where $r_n := \frac{n}{a_n}$.
\end{theorem}

\medskip

The rest of the paper is organized as follows:

\medskip


\qquad $\bullet$ Section \ref{s:preliminaries:BRW}, environment: preliminaries on branching random walks.


\qquad $\bullet$ Section \ref{s:proof-Zn}, environment: proof of Theorem \ref{t:Yn}.

\qquad $\bullet$ Section \ref{s:preliminaries:biased_walks}, biased walk: preliminaries on hitting barriers and local times.

\qquad $\bullet$ Section \ref{s:proof-p:x=0}, biased walk: proof  of Proposition \ref{p:x=0}.

\qquad $\bullet$ Section \ref{s:proof-t:ligne-arret}, biased walk: proof of Theorem \ref{t:ligne-arret}.

\qquad $\bullet$ Section \ref{s:localproba}, biased walk: proof  of Proposition  \ref{p:local}.

\qquad $\bullet$ Section \ref{s:maintheorem}, biased walk: proofs of Lemma \ref{l:pinpim}, Theorem \ref{t:cv} and Corollary \ref{c:cv}.

\medskip

Some comments on the organization are in order.  In the next two sections, we study the behaviour of the random environment, starting in Section \ref{s:preliminaries:BRW} by recalling some known results for branching random walks, and ending in Section \ref{s:proof-Zn} with the proof of Theorem \ref{t:Yn}. The biased walk $(X_n)$ comes into picture in the last five sections. In Section \ref{s:preliminaries:biased_walks}, we collect a couple of useful results about hitting lines and local times for the biased walk. The proof of Proposition \ref{p:x=0}, which is the most technical part of the paper, is presented in Section \ref{s:proof-p:x=0}. Once Proposition \ref{p:x=0} is established, we use it to deduce Theorem \ref{t:ligne-arret} in Section \ref{s:proof-t:ligne-arret}, and Proposition \ref{p:local} in Section \ref{s:localproba}. Finally, Theorem \ref{t:cv} and Corollary \ref{c:cv} (together with Lemma \ref{l:pinpim}) are proved in Section \ref{s:maintheorem}.

Throughout the paper, for any pair of vertices $x$ and $y$, we write $x<y$ or $y>x$ if $y$ is a (strict) descendant of $x$, and $x\le y$ or $y\ge x$ if either $y>x$, or $y=x$.

 \medskip

\section{Environment: preliminaries on branching random walks}
\label{s:preliminaries:BRW}


$\phantom{aob}$We recall, in this section, some known results in the literature for branching random walks, and deduce a few useful consequences. 

Under assumption (\ref{cond-hab}), there exists a sequence of i.i.d.\ real-valued random variables $(S_i-S_{i-1}, \, i\ge 0)$, with $S_0=0$, such that for any $n\ge 1$ and any Borel function $g: \r^n \to \r_+$,
\begin{equation}
    \E \Big[ \sum_{x\in \T: \, |x|=n} g(V(x_i), \, 1\le i\le n) \Big]
    =
    \E \Big[ \ee^{S_n} \, g(S_i, \, 1\le i\le n) \Big] \, ,
    \label{many-to-one}
\end{equation}

\noindent where, for any vertex $x\in \T$, $x_i$ ($0\le i\le n$) denotes, as before, the ancestor of $x$ in the $i$-th generation. As such, $V(x_0)$, $V(x_1)$, $\cdots$, $V(x_n)$ (for $|x|=n$) are the values of the potential $V$ alongs the branch $[\![ \varnothing , \, x]\!]$.

Formula (\ref{many-to-one}), often referred to as the ``many-to-one formula", is easily checked by induction on $n$. However, the appearance of the new, one-dimensional random walk $(S_i, \, i\ge 0)$ has a deep meaning in terms of the so-called spinal decomposition via a change of probabilities. The idea of change of probabilities in the study of spatial branching processes has a long history, going back at least to Kahane and Peyri\`ere~\cite{kahane-peyriere} and to Bingham and Doney~\cite{bingham-doney}, and has led to various forms of the spinal decomposition. Since Lyons, Pemantle and Peres~\cite{lyons-pemantle-peres}, it reaches a standard way of presentation. In our paper, we do not need any deep applications of the spinal decomposition, so we stay with the original probability $\P$ without making any change of probabilities, even though we do need a ``bivariate" version of (\ref{many-to-one}): 
\begin{equation}
    \E \Big[ \sum_{x\in \T: \, |x|=n} g\Big( V(x_i), \, \Lambda(x_{i-1}), \; 1\le i\le n\Big) \Big]
    =
    \E \Big[ \ee^{S_n} \, g\Big( S_i, \, \widetilde{\Lambda}_{i-1}, \; 1\le i\le n\Big) \Big] ,
    \label{spinal-decomposition}
\end{equation}

\noindent where, on the left-hand side of (\ref{spinal-decomposition}), we define
\begin{equation}
    \Lambda(x)
    :=
    \sum_{y: \, {\buildrel \leftarrow \over y} = x} \ee^{-[V(y)-V(x)]} \, ,
    \qquad
    x\in \T \, ,
    \label{Lambda}
\end{equation}

\noindent and on the right-hand side of (\ref{spinal-decomposition}), $(\widetilde{\Lambda}_i, \, i\ge 0)$ is such that $(S_i-S_{i-1}, \, \widetilde{\Lambda}_{i-1})$, $i\ge 1$, are i.i.d.\ random vectors whose law is characterized by
\begin{equation}
    \E \Big[ h ( S_1, \, \widetilde{\Lambda}_0) \Big]
    =
    \E \Big[ \sum_{x\in \T: \, |x|=1} \ee^{-V(x)} h(V(x), \, \Lambda(\varnothing) ) \Big] \, ,
    \label{joint-law:(S,Lambda)}
\end{equation}

\noindent for any Borel function $h: \, \r^2 \to \r_+$. The last equality follows, obviously, from (\ref{spinal-decomposition}) by taking $n=1$ there. Note that by definition, $\Lambda(\varnothing) = \sum_{x: \, |x|=1} \ee^{-V(x)}$. 

In particular, an application of the H\"older inequality, using assumption (\ref{integrability-assumption}), yields the existence of $\delta_1>0$ such that
\begin{equation}
    \E [ (\widetilde{\Lambda}_0)^{\delta_1} ]
    =
    \Big[ \Big( \sum_{x: \, |x|=1} \ee^{-V(x)} \Big)^{1+\delta_1} \Big]
    <
    \infty\, .
    \label{integrabilite_sous_Q}
\end{equation}

These are known facts about the spinal decomposition. For a proof of (\ref{spinal-decomposition}), see \cite{yzenergy}. 

We now deduce several simple but useful results. The first allows us to include the random variable $\Lambda(x)$ in the bivariate many-to-one formula (\ref{spinal-decomposition}). The second takes care of summation over all vertices on the stopping line $\mathscr{L}_r$ instead of on a given generation, which leads to the third which is also the main estimate in this section.

\medskip

\begin{lemma}
\label{l:many-to-one-appli1}

 Assume $(\ref{cond-hab})$. 
 Let $\Lambda(x)$ be as in $(\ref{Lambda})$.
 For any $n\ge 1$ and any Borel function $f: \, \r^{2n+1} \to \r_+$,
 we have
 $$
 \E \Big[ \sum_{x\in \T: \, |x|=n} f\Big( V(x_i), \, \Lambda(x_{i-1}), \; 1\le i\le n, \, \Lambda(x) \Big) \Big]
 =
 \E \Big[ \ee^{S_n} \, F\Big( S_i, \, \widetilde{\Lambda}_{i-1}, \; 1\le i\le n\Big) \Big] \, ,
 $$
 where $(S_i-S_{i-1}, \, \widetilde{\Lambda}_{i-1})$, $i\ge 1$, are i.i.d.\ whose common distribution is given in $(\ref{joint-law:(S,Lambda)})$, and
 \begin{equation}
     F(a_i, \, b_{i-1}, \; 1\le i\le n)
     :=
     \E \Big[ f\Big(a_i, \, b_{i-1}, \; 1\le i\le n, \, \sum_{x\in \T: \, |x|=1} \ee^{-V(x)}\Big) \Big] \, .
     \label{F}
 \end{equation}
 
 \noindent In particular, if $g: \r^{n+1} \to \r_+$ is a Borel function, then
 $$
 \E \Big[ \sum_{x\in \T: \, |x|=n} g\Big( V(x_1), \, \cdots, \, V(x_n), \, \Lambda(x) \Big) \Big]
 =
 \E \Big[ \ee^{S_n} \, G\Big( S_1, \, \cdots, \, S_n\Big) \Big] \, ,
 $$
 where $G(a_1, \, \cdots, \, a_n) := \E  [ g(a_1, \, \cdots, \, a_n, \, \sum_{x\in \T: \, |x|=1} \ee^{-V(x)}) ]$.
  
\end{lemma}

\medskip

\noindent {\it Proof.} Let $\mathscr{F}_n := \sigma(x, \, V(x), \, x\in \T, \, |x| \le n)$, the $\sigma$-field generated by the branching random walk in the first $n$ generations. By definition, for $|x|=n$, $\Lambda(x)$ is independent of $\mathscr{F}_n$, so
\begin{eqnarray*}
 &&\E \Big[ \sum_{x\in \T: \, |x|=n} f\Big( V(x_i), \, \Lambda(x_{i-1}), \; 1\le i\le n, \, \Lambda(x) \Big) \, \Big| \, \mathscr{F}_n \Big]
    \\
 &=&\sum_{x\in \T: \, |x|=n} F\Big( V(x_i), \, \Lambda(x_{i-1}), \; 1\le i\le n\Big) \, ,
\end{eqnarray*}

\noindent where $F$ is given by (\ref{F}). Taking expectation with respect to $\P$ on both sides, and using the bivariate many-to-one formula (\ref{spinal-decomposition}), we obtain the lemma.\hfill$\Box$

\medskip

\begin{lemma}
\label{l:many-to-one-appli2}

 Assume $(\ref{cond-hab})$. 
 Let $E_1$, $E_2$, $\cdots$ be Borel subsets of $\r$.
 Let $r>1$ and let $\mathscr{L}_r$ be as in $(\ref{gamma})$.
 Then
 \begin{equation}
     \E\Big( \sum_{x\in \mathscr{L}_r} \ee^{-V(x)} \, {\bf 1}_{\{ V(x_i) \in E_i , \, 1\le i \le |x|\} } \Big)
     =
     \P \Big( S_i \in E_i , \, 1\le i \le T_r^{(S)} \Big) \, ,
     \label{many-to-one-appli-1}
 \end{equation}

 \noindent where $T_r^{(S)} := \inf\{ i \ge 1: \, \sum_{j=1}^i \ee^{S_j-S_i} > r\}$. 

\end{lemma}

\medskip

\noindent {\it Proof.} We write
$$
\sum_{x\in \mathscr{L}_r} \ee^{-V(x)} \, {\bf 1}_{\{ V(x_i) \in E_i , \, 1\le i \le |x|\} }
=
\sum_{k=1}^\infty \sum_{x: \, |x|=k} \ee^{-V(x)} \, {\bf 1}_{\{ x\in \mathscr{L}_r\} } \, {\bf 1}_{\{ V(x_i) \in E_i , \, 1\le i \le k\} } \, .
$$

\noindent Obviously, $\{ x\in \mathscr{L}_r\} = \{ \sum_{z\in \, ]\! ] \varnothing, \, x]\! ]} \ee^{V(z)-V(x)} > r , \; \sum_{z\in \, ]\! ] \varnothing, \, v]\! ]} \ee^{V(z)-V(v)} \le r, \, \forall v\in \, ]\!] \varnothing, \, x[\![ \, \}$. We take expectation with respect to $\P$ on both sides. By the many-to-one formula (\ref{many-to-one}), 
$$
\E\Big( \sum_{x\in \mathscr{L}_r} \ee^{-V(x)} \, {\bf 1}_{\{ V(x_i) \in E_i , \, 1\le i \le |x|\} } \Big)
=
\sum_{k=1}^\infty \P\Big( T_r^{(S)} = k, \; S_i \in E_i , \, 1\le i \le k \Big) \, ,
$$

\noindent which is $\P ( S_i \in E_i , \, 1\le i \le T_r^{(S)})$.\hfill$\Box$

\bigskip

\begin{remark}
\label{r:E(Yn)}

Let $Y_r :=\sum_{x\le \mathscr{L}_r}\ee^{-V(x)} = 1+ \sum_{k=1}^\infty \sum_{x: \, |x|=k} \ee^{-V(x)} \, {\bf 1}_{\{ x\le \mathscr{L}_r\} }$ as in \eqref{Yn}. Since $x\le \mathscr{L}_r$ means $\sum_{z\in \, ]\! ] \varnothing, \, v]\! ]} \ee^{V(z)-V(v)} \le r$ for all $v\in \, ]\!] \varnothing, \, {\buildrel \leftarrow \over x}[\![ \,$ (the inequality considered as holding trivially if $|x| =1$), the proof of Lemma \ref{l:many-to-one-appli2} yields 
$$
\E(Y_r)
=
1+ \E (T_r^{(S)})
\le 1+  \E \Big[ \inf\{ i\ge 1: \, \max_{1\le j\le i}S_j - S_i > \log r \} \Big]  \, .
$$

\noindent It is easy to check (for a detailed proof, see \cite{yzenergy}) that $\E [\, \inf\{ i\ge 1: \, \max_{1\le j\le i}S_j - S_i > u\}\, ]$ is bounded by $c_1 \, u^2$ for some constant $c_1>0$ and all $u\ge 1$. Hence, there exists $c_2>0$ such that
\begin{equation}
    \E(Y_r) 
    \le 
    c_2 \, (\log r)^2 \, ,
    \qquad r\ge 2\, .
    \label{E(Yn)<}
\end{equation}

\noindent We are going to use (\ref{E(Yn)<}) in Section \ref{s:proof-p:x=0}, in the proof of Proposition \ref{p:x=0}.

Although we do not need it in the present paper, an elementary argument shows that $\frac{\E(Y_r)}{(\log r)^2}$ is bounded also from below.
\hfill$\Box$

\end{remark}

\bigskip

We now present the main probabilistic estimate of the section.

\medskip

\begin{lemma}
\label{l:sum(xinLn)}

 Assume $(\ref{cond-hab})$ and $(\ref{integrability-assumption})$. The laws of $(\log r) \sum_{x\in \mathscr{L}_r}\ee^{-V(x)}$ under $\P^*$, for $r\ge 1$, are tight.

 In particular, for any $a<1$, $(\log r)^a \sum_{x\in \mathscr{L}_r}\ee^{-V(x)} \to 0$, $r\to \infty$, in $\P^*$-probability.

\end{lemma}

\medskip

\noindent {\it Proof.} Let $\varepsilon>0$. Our assumption ensures $\inf_{x: \, |x| = n} V(x) \to \infty$ (for $n\to \infty$) $\P^*$-a.s.\ (see (\ref{U->infty}); so we can choose and fix a constant $\alpha >0$ such that
\begin{equation}
    \P^* \Big( \inf_{x\in \T} V(x) \ge -\alpha \Big)
    \ge
    1-\varepsilon .
    \label{alpha}
\end{equation}

\noindent For any $x\in \T$, write
$$
\underline{V}(x)
:=
\min_{y\in [\![ \varnothing, \, x]\!]} V(y) \, .
$$

\noindent By Lemma \ref{l:many-to-one-appli2},
\begin{equation}
    \E\Big( 
    \sum_{x\in \mathscr{L}_r}
    \ee^{-V(x)} \, 
    {\bf 1}_{\{ \underline{V}(x) \ge -\alpha\} } \Big)
    =
    \P \Big( \underline{S}_{T_r^{(S)}} \ge -\alpha \Big),
    \label{pf-partition-fct}
\end{equation}

\noindent where $T_r^{(S)} := \inf\{ i \ge 1: \, \sum_{j=1}^i \ee^{S_j-S_i} > r\}$, and $\underline{S}_i := \min_{0\le j\le i} S_j$. 

Let $H(\frac12 \log r) := \inf\{ i\ge 1: \, S_i > \frac12 \log r\}$. We have
$$
\P \Big( \underline{S}_{T_r^{(S)}} \ge -\alpha \Big)
\le
\P \Big( T_r^{(S)} < H(\frac12 \log r), \; \underline{S}_{T_r^{(S)}} \ge -\alpha \Big)
+
\P \Big( \underline{S}_{H(\frac12 \log r)} \ge -\alpha \Big) \, .
$$

We bound the two probability expressions on the right-hand side. For $\P ( \underline{S}_{H(\frac12 \log r)} \ge -\alpha)$, we write $H_-(-\alpha) := \inf \{ i\ge 1: \,  S_i < -\alpha\}$, to see that for some constant $c_3>0$,
$$
\P \Big( \underline{S}_{H(\frac12 \log r)} \ge -\alpha \Big)
=
\P \{ H(\frac12 \log r) < H_-(-\alpha) \}
\le
\frac{c_3 \, \alpha}{  \frac12 \log r + \alpha} \, .
$$

\noindent [For the last inequality, which is elementary, see for example A\"\i d\'ekon~\cite{aidekon2010ECP} under the assumption of existence of exponential moments of $S_1$.] Hence $\P ( \underline{S}_{H(\frac12 \log r)} \ge -\alpha) \le \frac{2 c_3 \, \alpha}{\log r}$. Accordingly,
\begin{equation}
    \P \Big( \underline{S}_{T_r^{(S)}} \ge -\alpha \Big)
    \le
    \P \Big( T_r^{(S)} < H(\frac12 \log r), \; \underline{S}_{T_r^{(S)}} \ge -\alpha \Big)
    +
   \frac{2 c_3 \, \alpha}{\log r} \, .
    \label{alpha1}
\end{equation}

\noindent To deal with $\P ( T_r^{(S)} < H(\frac12 \log r), \; \underline{S}_{T_r^{(S)}} \ge -\alpha )$, we note that by definition of $T_r^{(S)}$, $r < \sum_{j=1}^{T_r^{(S)}} \ee^{S_j-S_{T_r^{(S)}}}$, which, on the event $\{ T_r^{(S)} < H(\frac12 \log r), \; \underline{S}_{T_r^{(S)}} \ge -\alpha\}$, is bounded by $\sum_{j=1}^{T_r^{(S)}} \ee^{S_j+\alpha} \le \sum_{j=1}^{H(\frac12 \log r) -1} \ee^{(\frac12 \log r) +\alpha} \le r^{1/2}\, \ee^\alpha \, H(\frac12 \log r)$. Consequently,
$$
\P \Big( T_r^{(S)} < H(\frac12 \log r), \; \underline{S}_{T_r^{(S)}} \ge -\alpha \Big)
\le
\P\Big( H(\frac12 \log r) > r^{1/2}\, \ee^{-\alpha} \Big) \, .
$$

\noindent By Kozlov~\cite{kozlov}, $\P\{ H(\frac12 \log r) > r^{1/2}\, \ee^{-\alpha} \} \le c_4 \frac{\ee^{\alpha/2} \log r}{r^{1/4}}$ for some constant $c_4>0$ and all $n\ge 2$. Going back to (\ref{alpha1}) and having (\ref{pf-partition-fct}) in mind, we obtain:
$$
\E\Big( \sum_{x\in \mathscr{L}_r}\ee^{-V(x)} {\bf 1}_{\{ \underline{V}(x) \ge -\alpha\} }\Big)
\le
c_4 \frac{\ee^{\alpha/2} \log r}{r^{1/4}}
+
\frac{2 c_3 \, \alpha}{\log r} \, .
$$

\noindent In view of (\ref{alpha}), and since $\varepsilon>0$ can be as small as possible, Lemma \ref{l:sum(xinLn)} follows readily.\hfill$\Box$

\section{Environment: proof of Theorem \ref{t:Yn}}
\label{s:proof-Zn}

$\phantom{aob}$This section is mainly devoted to the proof of Theorem \ref{t:Yn}, but also prepares a few useful estimates for the forthcoming sections. The material in this section concerns only the environment (thus the potential $V$ and the symmetrized potential $U$); no discussion on the movement of the biased walk $(X_i)$ is involved.

Let $W_n := \sum_{|x|=n} \ee^{-V(x)}$, $n\ge 0$, be the additive martingale as in (\ref{Wn}). Consider also, for $n\ge 0$ and $\lambda >0$,
\begin{equation}
    W_n^{(\lambda)}
    :=
    \sum_{|x|=n} \ee^{-V(x)}\, 
    {\bf 1}_{ \{ \max_{y\in \, ]\! ]\varnothing, \, x]\! ]} \, [\overline{V}(y) - V(y)] \le \lambda \} } \, ,
    \label{Wn(alpha)}
\end{equation}

\noindent where
$$
\overline{V}(y) 
:=
\max_{z\in [ \! [ \varnothing, \, y]\! ]} V(z) \, .
$$

\noindent We mentioned earlier in \eqref{W->0} that under assumption (\ref{cond-hab}), we have $W_n \to 0$ $\P^*$-a.s. The rate of decay of $W_n$ is known: according to \cite{ezratio}, under (\ref{cond-hab}) and (\ref{integrability-assumption}), we have
\begin{equation}
    \lim_{n\to \infty} \, n^{1/2} \, W_n
    =
    \Big( \frac{2}{\pi \sigma^2} \Big)^{1/2} D_\infty \, ,
    \qquad \hbox{in $\P^*$-probability.}
    \label{ezratio}
\end{equation}
 
\noindent The asymptotics of $W_n^{(\lambda)}$ are also studied: according to Madaule~\cite{madaule}, for any $a\ge 0$,
$$
    \lim_{n\to \infty} \, \frac{W_n^{(n^{1/2}a)}}{W_n}
    =
    \P \Big( \eta < \frac{a}{\sigma} \Big) \, ,
    \qquad \hbox{in $\P^*$-probability,}
$$

\noindent where $\eta := \sup_{s\in [0, \, 1]} [\, \overline{\mathfrak{m}}(s) - \mathfrak{m}(s)]$, with $\overline{\mathfrak{m}}(s) := \sup_{u\in [0, \, s]} \mathfrak{m}(u)$, and $(\mathfrak{m}(s), \, s\in [0, \, 1])$ denoting as before a standard Brownian meander. In view of (\ref{ezratio}), this is equivalent to saying the following convergence in $\P^*$-probability:
\begin{equation}
    \lim_{n\to \infty} \, n^{1/2} \, W_n^{(n^{1/2}a)} 
    =
    \Big( \frac{2}{\pi \sigma^2} \Big)^{1/2} D_\infty \,
    \P \Big( \eta < \frac{a}{\sigma} \Big) .
    \label{madaule}
\end{equation}

\noindent This holds for any given $a\ge 0$. 

By the absolute continuation relation between the Brownian meander $\mathfrak{m}$ and the three-dimensional Bessel process $R$ (see   \cite{imhof}), we know that the law of $\eta := \sup_{s\in [0, \, 1]} [\, \overline{\mathfrak{m}}(s) - \mathfrak{m}(s)]$ is absolutely continuous with respect to the law of $\sup_{s\in [0, \, 1]} [\, \overline{R}(s) - R(s)]$. The latter is atomless because $R$ is an $h$-transform (in the sense of Doob) of Brownian motion. As a consequence, $a\mapsto \P ( \eta < \frac{a}{\sigma})$ is continuous on $\r$. On the other hand, both $a\mapsto W_n^{(n^{1/2}a)}$ and $a\mapsto \P ( \eta < \frac{a}{\sigma})$ are non-decreasing. It follows that (\ref{madaule}) holds uniformly (in $a\ge 0$) in the following sense: for any $\varepsilon>0$,
\begin{equation}
    \lim_{n\to \infty} \, 
    \P^* \Big\{ \sup_{a\ge 0} 
    \Big| 
    n^{1/2} \, W_n^{(n^{1/2}a)} 
    -
    \Big( \frac{2}{\pi \sigma^2} \Big)^{1/2} D_\infty \, \P ( \eta < \frac{a}{\sigma})
    \Big|
    \ge \varepsilon \Big\}
    =
    0 \, .
    \label{madaule2}
\end{equation}

\noindent We now state a lemma.

\medskip

\begin{lemma}
\label{l:sum_martingales}

 We have
 \begin{equation}
     \lim_{\lambda\to \infty} \, \frac{1}{\lambda} \sum_{k=1}^\infty W_k^{(\lambda)}
     =
     \frac{2}{\sigma^2} \, D_\infty \, ,
     \qquad\hbox{in $\P^*$-probability} \, .
     \label{Zn-proof1}
 \end{equation}
\end{lemma}

\medskip

\noindent {\it Proof.} We first argue that in $\sum_{k=1}^\infty W_k^{(\lambda)}$, only those $k$ that are comparable to $\lambda^2$ make a significant contribution to the sum. More precisely, we claim that for any $\varepsilon_1>0$, 
\begin{eqnarray}
    \lim_{b\to 0} \, \limsup_{\lambda\to \infty} \, \P^* \Big\{ \frac{1}{\lambda} \sum_{k=1}^{\lfloor b\lambda^2\rfloor} W_k^{(\lambda)} \ge \varepsilon_1 \Big\}
 &=& 0 \, ,
    \label{Zn-proof2}
    \\
    \lim_{B\to \infty} \, \limsup_{\lambda\to \infty} \, \P^* \Big\{ \frac{1}{\lambda} \sum_{k= \lfloor B\lambda^2\rfloor}^\infty W_k^{(\lambda)} \ge \varepsilon_1 \Big\}
 &=& 0 \, .
    \label{Zn-proof3}    
\end{eqnarray}

To prove (\ref{Zn-proof2}) and (\ref{Zn-proof3}), let $\varepsilon>0$ and fix $\alpha\ge 0$ as in (\ref{alpha}), i.e., such that 
$$
\P^*\{ \inf_{x\in \T} V(x) \ge -\alpha\} 
\ge 
1-\varepsilon \, .
\leqno(\ref{alpha})
$$

\noindent Consider the truncated version of $W_k^{(\lambda)}$ defined by
$$
W_k^{(\lambda,\, \alpha)}
:=
\sum_{|x|=k} \ee^{-V(x)}\, 
{\bf 1}_{ \{ \max_{y\in \, ]\! ] \varnothing, \, x]\! ]} \, [\overline{V}(y) - V(y)] \le \lambda \} } \,
{\bf 1}_{\{ \underline{V}(x) \ge -\alpha\} } \, ,
$$

\noindent where $\underline{V}(x) := \min_{z\in [\![ \varnothing, \, x]\!]} V(z)$ as before. 
 Clearly,  on the set $\{ \inf_{x\in \T} V(x) \ge -\alpha\}$, $W_k^{(\lambda,\, \alpha)}= W_k^{(\lambda)}$  for all $k \ge1$.

By the many-to-one formula in (\ref{many-to-one}), 
\begin{equation}
    \E(W_k^{(\lambda,\, \alpha)})
    =
    \P \Big\{ \max_{0\le j\le k} (\overline{S}_j - S_j) \le \lambda, \, \underline{S}_k \ge -\alpha \Big\} \, ,
    \label{Zn-proof4}
\end{equation}

\noindent where $\overline{S}_j := \max_{0\le i\le j} S_i$ and $\underline{S}_j := \min_{0\le i\le j} S_i$.

The proof of (\ref{Zn-proof2}) is easy: we have $\E(W_k^{(\lambda,\, \alpha)}) \le \P \{ \underline{S}_k \ge -\alpha\}$, which is bounded by $\frac{c_5}{k^{1/2}}$ for some constant $c_5$ (depending on $\alpha$) and all $k\ge 1$ (see Kozlov~\cite{kozlov}); hence
$$
\frac{1}{\lambda} \E \Big( \sum_{k=1}^{\lfloor b\lambda^2\rfloor} W_k^{(\lambda, \, \alpha)} \Big)
\le
\frac{1}{\lambda} \sum_{k=1}^{\lfloor b\lambda^2\rfloor} \frac{c_5}{k^{1/2}} \, ,
$$

\noindent which goes to $0$ when $\lambda \to \infty$ and then $b\to 0$. This readily yields (\ref{Zn-proof2}).

To prove (\ref{Zn-proof3}), we use (\ref{Zn-proof4}), and apply the Markov property at time $\frac{k}{2}$ (treating it as an integer by dropping the symbol of the integer part), to see that
\begin{eqnarray*}
    \E(W_k^{(\lambda,\, \alpha)})
 &\le&\P \Big\{ \underline{S}_{\frac{k}{2}} \ge -\alpha, \; \max_{\frac{k}{2}\le j\le k} (\overline{S}_j - S_j) \le \lambda \Big\} 
    \\
 &\le& \P \{ \underline{S}_{\frac{k}{2}} \ge -\alpha \} 
    \times
    \P \Big\{ \max_{0\le j\le \frac{k}{2}} (\overline{S}_j - S_j) \le \lambda \Big\} \, .
\end{eqnarray*}

\noindent Again, $\P \{ \underline{S}_{\frac{k}{2}} \ge -\alpha \} \le \frac{c_5}{(\frac{k}{2})^{1/2}}$, whereas $\P \{ \max_{0\le j\le \frac{k}{2}} (\overline{S}_j - S_j) \le \lambda \}$ can be estimated as follows: by the Markov property, $\P \{ \max_{0\le j\le \frac{k}{2}} (\overline{S}_j - S_j) \le \lambda \} \le [ \, \P \{ \max_{0\le j\le \lfloor \lambda^2\rfloor} (\overline{S}_j - S_j) \le \lambda \} \, ]^{\lfloor \frac{k}{2\, \lfloor \lambda^2\rfloor} \rfloor}$. By Donsker's theorem, there exists a constant $0<c_6<1$ such that $\P \{ \max_{0\le j\le \lfloor \lambda^2\rfloor} (\overline{S}_j - S_j) \le \lambda \} \le 1-c_6$ for all sufficiently large $\lambda$ (say $\lambda\ge \lambda_0$) which yields $\P \{ \max_{0\le j\le \frac{k}{2}} (\overline{S}_j - S_j) \le \lambda \} \le (1-c_6)^{\lfloor \frac{k}{2\, \lfloor \lambda^2\rfloor} \rfloor}$, $\forall \lambda\ge \lambda_0$. Hence, for $\lambda \ge \lambda_0$ and $k\ge 1$,
$$
\E(W_k^{(\lambda, \, \alpha)})
\le
\frac{c_5}{(\frac{k}{2})^{1/2}} \, (1-c_6)^{\lfloor \frac{k}{2\, \lfloor \lambda^2\rfloor} \rfloor} \, ,
$$

\noindent from which it follows that
$$
\lim_{B\to \infty} \, \limsup_{\lambda\to \infty} \, \frac{1}{\lambda} \, \E \Big\{ \sum_{k= \lfloor B\lambda^2\rfloor}^\infty W_k^{(\lambda, \, \alpha)} \Big\} 
= 0 \, .
$$

\noindent Since $\varepsilon>0$ in (\ref{alpha}) can be as small as possible, this implies (\ref{Zn-proof3}).

Now that (\ref{Zn-proof2}) and (\ref{Zn-proof3}) are justified, we are ready for the proof of Lemma \ref{l:sum_martingales}. Fix $B>b>0$. By (\ref{madaule2}), for $\lambda\to \infty$,
$$
    \frac{1}{\lambda} \sum_{k= \lfloor b\lambda^2\rfloor}^{\lfloor B\lambda^2\rfloor} W_k^{(\lambda)}
    =
    \frac{1}{\lambda} \Big( \frac{2}{\pi \sigma^2} \Big)^{1/2} D_\infty \sum_{k= \lfloor b\lambda^2\rfloor}^{\lfloor B\lambda^2\rfloor} \frac{1}{k^{1/2}}\, \P \Big( \eta < \frac{\lambda}{\sigma\, k^{1/2}} \Big) 
    +
    o_{\P^*}(1) \, ,
$$

\noindent where $o_{\P^*}(1)$ denotes a term satisfying $\lim_{\lambda\to \infty} o_{\P^*}(1) =0$ in $\P^*$-probability (whose value may vary from line to line). On the other hand, by Fubini's theorem,
\begin{eqnarray*}
    \frac{1}{\lambda} \int_{b\lambda^2}^{B\lambda^2} \frac{1}{u^{1/2}}\, \P \Big( \eta < \frac{\lambda}{\sigma\, u^{1/2}} \Big) \d u
 &=& \frac{1}{\lambda} \, \E\Big[ \int_{b\lambda^2}^{(B\lambda^2) \wedge \frac{\lambda^2}{\sigma^2 \eta^2}}  \frac{\d u }{u^{1/2}}\, {\bf 1}_{\{ \eta < \frac{1}{\sigma\, b^{1/2}} \} }  \Big]
    \\
 &=&2 \, \E\Big[ \Big( (B^{1/2} \wedge \frac{1}{\sigma \eta} ) - b^{1/2} \Big) \, {\bf 1}_{\{ \eta < \frac{1}{\sigma\, b^{1/2}} \} } \Big] \, .
\end{eqnarray*}

\noindent Since $\eta$ is atomless, this yields
\begin{equation}
    \frac{1}{\lambda} \sum_{k= \lfloor b\lambda^2\rfloor}^{\lfloor B\lambda^2\rfloor} W_k^{(\lambda)}
    =
    \Big( \frac{8}{\pi \sigma^2} \Big)^{1/2} D_\infty\, \E\Big[ \Big( (B^{1/2} \wedge \frac{1}{\sigma \eta} ) - b^{1/2} \Big) \, {\bf 1}_{\{ \eta < \frac{1}{\sigma\, b^{1/2}} \} } \Big] 
    +
    o_{\P^*}(1) \, .
    \label{Zn-proof5}
\end{equation}

\noindent Note that $\E [ ( (B^{1/2} \wedge \frac{1}{\sigma \eta} ) - b^{1/2} ) \, {\bf 1}_{\{ \eta < \frac{1}{\sigma\, b^{1/2}} \} } ]  \to \E [\frac{1}{\sigma \eta}]$ when $B\to \infty$ and $b\to 0$. In view of (\ref{Zn-proof2}) and (\ref{Zn-proof3}), we see that when $\lambda\to \infty$, 
$$
\frac{1}{\lambda} \sum_{k=1}^\infty W_k^{(\lambda)}
\to
\Big( \frac{8}{\pi \sigma^2} \Big)^{1/2} D_\infty \, \E [\frac{1}{\sigma \eta}] \, ,
\qquad\hbox{in $\P^*$-probability} \, .
$$

\noindent By \cite{Yor}, $\E (\frac{1}{\eta}) = (\frac{\pi}{2})^{1/2}$, which yields Lemma \ref{l:sum_martingales}.\hfill$\Box$

\bigskip

We now have all the ingredients for the proof of Theorem \ref{t:Yn}.

\bigskip

\noindent {\it Proof of Theorem \ref{t:Yn}.} By definition,
$$
Y_r
=
\sum_{x\in \T} \ee^{-V(x)} \, {\bf 1}_{\{ x< \mathscr{L}_r\} }
+
\sum_{x\in \mathscr{L}_r} \ee^{-V(x)} .
$$

\noindent We already know (Lemma \ref{l:sum(xinLn)}) that $\sum_{x\in \mathscr{L}_r} \ee^{-V(x)} \to 0$ in $\P^*$-probability. So it remains to check that
\begin{equation}
    \lim_{r\to \infty} \, \frac{1}{\log r} \sum_{x\in \T} \ee^{-V(x)} \, {\bf 1}_{\{ x<\mathscr{L}_r\} }
    =
    \frac{2}{\sigma^2} \, D_\infty \, ,
    \qquad\hbox{in $\P^*$-probability} \, .
    \label{Zn-step1-proof-2}
\end{equation}

By definition, $\{x<\mathscr{L}_r\}$ means $\sum_{z\in \, ]\! ]\varnothing, \, y]\! ]} \ee^{V(z)-V(y)} \le r$, $\forall y\in \, ]\!] \varnothing, \, x]\!]$. So
\begin{eqnarray}
    \sum_{x\in \T} \ee^{-V(x)} \, {\bf 1}_{\{ x< \mathscr{L}_r\} }
 &=&\sum_{k=0}^\infty \sum_{x: \, |x|=k} \ee^{-V(x)} \, {\bf 1}_{\{ \sum_{z\in \, ]\! ]\varnothing, \, y]\! ]} \ee^{V(z)-V(y)} \le r, \; \forall y\in \, ]\!] \varnothing, \, x]\!]\} }
    \nonumber
    \\
 &\le& \sum_{k=0}^\infty \sum_{x: \, |x|=k} \ee^{-V(x)} \, {\bf 1}_{\{ \max_{y\in \, ]\! ]\varnothing, \, x]\! ]} \, [\overline{V}(y) - V(y)] \le \log r\} }
    \nonumber
    \\
 &=&\sum_{k=0}^\infty W_k^{(\log r)} \, .
    \label{pf-Yn-ub}
\end{eqnarray}

A similar lower bound holds as well: we fix an arbitrary positive real number $B>0$,
\begin{eqnarray}
    \sum_{x\in \T} \ee^{-V(x)} \, {\bf 1}_{\{ x< \mathscr{L}_r\} }
 &\ge& \sum_{k=0}^{\lfloor B(\log r)^2 \rfloor} \sum_{x: \, |x|=k} \ee^{-V(x)} \, {\bf 1}_{\{ \max_{y\in \, ]\! ] \varnothing, \, x]\! ]} \, [\overline{V}(y) - V(y)] \le \log\frac{r}{B(\log r)^2} \} } 
    \nonumber
    \\
 &=&\sum_{k=0}^{\lfloor B(\log r)^2 \rfloor} W_k^{(\log\frac{r}{B (\log r)^2})} \, . \label{pf-Yn-low}
\end{eqnarray}

\noindent Applying Lemma \ref{l:sum_martingales} and (\ref{Zn-proof3}), and noting that $\lim_{r\to \infty} \frac{\log\frac{r}{B (\log r)^2}}{\log r} = 1$, we obtain that under $\P^*$,
$$
\lim_{r\to \infty} \, \frac{1}{\log r} \sum_{x\in \T} \ee^{-V(x)} \, {\bf 1}_{\{ x<\mathscr{L}_r\} }
=
\frac{2}{\sigma^2} \, D_\infty \, ,
\qquad\hbox{in probability} \, .
$$

\noindent Theorem \ref{t:Yn} is proved.\hfill$\Box$

\medskip

\begin{remark}

The proof of the upper bound in Theorem \ref{t:Yn}, combined with (\ref{Zn-proof3}), tells us that for any $\varepsilon>0$,
\begin{equation}
    \lim_{B\to \infty} \, \limsup_{r\to \infty} \, 
    \P^* \Big\{ \frac{1}{\log r} \sum_{x\in \T: \, |x| \ge B(\log r)^2, \, x\le \mathscr{L}_r} \ee^{-V(x)} \ge \varepsilon\Big\}
    =
    0 \, .
    \label{Zn-proof6}
\end{equation}

\noindent We are entitled to sum over $x\le \mathscr{L}_r$ instead of over $x<\mathscr{L}_r$ because $\sum_{x\in \mathscr{L}_r} \ee^{-V(x)} \to 0$ in $\P^*$-probability (Lemma \ref{l:sum(xinLn)}); (\ref{Zn-proof6}) will be useful in Section \ref{s:proof-p:x=0}.\hfill$\Box$  

\end{remark}

\section{Biased walks: preliminaries on hitting barriers and local times}
\label{s:preliminaries:biased_walks}

$\phantom{aob}$In this section, we collect two preliminary results for the biased walk $(X_i)$. The first is a weaker version of Theorem \ref{t:ligne-arret}, and the second concerns the covariance of edge local times. For the sake of clarity, we present them is two distinct subsections.

\subsection{Hitting reflecting barriers}
\label{subs:stopping_lines}

$\phantom{aob}$This subsection is devoted to a weaker version of Theorem \ref{t:ligne-arret}, stated as follows. The proof of Theorem \ref{t:ligne-arret} comes much later, in Section \ref{s:proof-t:ligne-arret}.

\begin{lemma}
\label{l:ligne-arret}

 Assume $(\ref{cond-hab})$ and $(\ref{integrability-assumption})$. If $r=r(n) := \frac{n}{(\log n)^\gamma}$ with $\gamma<1$, then
 $$
 \lim_{n\to \infty} \, 
 \p \Big( \bigcup_{i=1}^n\{ X_i \in \mathscr{L}_r \} \Big)
 =
 0 \, ,
 $$
 where $\mathscr{L}_r$ is as in $(\ref{gamma})$.

\end{lemma}

\medskip

\noindent {\it Proof.} Define
\begin{eqnarray}
    T_x
 &:=& \inf\{ i\ge 0: \, X_i =x \} \, , 
    \qquad x\in \T \, ,
    \label{T}
    \\
    T_\varnothing^+
 &:=& \inf\{ i\ge 1: \, X_i = \varnothing \} \, .
    \label{tau}
\end{eqnarray}

\noindent In words, $T_x$ is the first hitting time at $x$ by the biased walk, whereas $T_\varnothing^+$ is the first {\it return} time to the root $\varnothing$. 

Let $x\in \T \backslash\{ \varnothing\}$. The probability $P_\omega (T_x < T_\varnothing^+)$ only involves a one-dimensional random walk in random environment (namely, the restriction at $[\! [ \varnothing, \, x[\![\,$ of the biased walk $(X_i)$), so a standard result for one-dimensional random walks in random environment (Golosov~\cite{golosov}) tells us that\begin{equation}
    P_\omega (T_x < T_\varnothing^+)
    =
    \frac{\omega(\varnothing, \, x_1)\, \ee^{V(x_1)}}{\sum_{z\in \, ]\!] \varnothing, \, x]\!]} \ee^{V(z)}}    =
    \frac{\omega(\varnothing, \, {\buildrel \leftarrow \over \varnothing})}{\sum_{z\in \, ]\!] \varnothing, \, x]\!]} \ee^{V(z)}} \, ,
    \label{1D-MAMA}
\end{equation}

\noindent where $x_1$ is the ancestor of $x$ in the first generation.

Define $T_\varnothing^{(0)} := 0$ and inductively $T_\varnothing^{(k)} := \inf\{ i> T_\varnothing^{(k-1)}: \, X_i = \varnothing \}$, $k\ge 1$. In words, $T_\varnothing^{(k)}$ is the $k$-th return time of the biased walk $(X_i)$ to the root $\varnothing$. [In particular, $T_\varnothing^{(1)} = T_\varnothing^+$.] For $n\ge 1$, we have
\begin{eqnarray*}
    P_\omega (T_x \le n)
 &=&\sum_{k=0}^\infty P_\omega \Big[ T_x \le n, \; T_\varnothing^{(k)} \le T_x < T_\varnothing^{(k+1)} \Big]
    \\
 &=&\sum_{k=0}^\infty E_\omega \Big[ {\bf 1}_{\{ T_\varnothing^{(k)} \le n\} } \, P_\omega (T_x < T_\varnothing, \; T_x \le n-j)\Big|_{j:= T_\varnothing^{(k)}}\Big]
    \\
 &\le&\sum_{k=0}^\infty E_\omega \Big[ {\bf 1}_{\{ T_\varnothing^{(k)} \le n\} } \, P_\omega (T_x < T_\varnothing^+)\Big]
    \\
 &=& P_\omega (T_x < T_\varnothing^+) \, E_\omega (L_n(\varnothing)+1),
\end{eqnarray*}

\noindent where $L_n (\varnothing) := \sum_{i=1}^n {\bf 1}_{\{ X_i = \varnothing \} }$ is the local time at $\varnothing$. By (\ref{1D-MAMA}), we get
$$
    P_\omega (T_x \le n)
    \le
    \frac{  E_\omega (L_n(\varnothing)+1)}{\sum_{z\in \, ]\! ] \varnothing, \, x]\! ]} \ee^{V(z)}} .
$$

Let $r>1$, and let $\mathscr{L}_r$ be as in (\ref{gamma}). We have
$$
P_\omega \Big( \bigcup_{i=1}^n\{ X_i \in \mathscr{L}_r \} \Big)
\le
\sum_{x\in \mathscr{L}_r} P_\omega (T_x \le n)
\le
E_\omega (L_n(\varnothing)+1) \sum_{x\in \mathscr{L}_r} \frac{1}{\sum_{z\in \, ]\! ] \varnothing, \, x]\! ]} \ee^{V(z)}} .
$$

\noindent By definition of $\mathscr{L}_r$, $\frac{1}{\sum_{z\in \, ]\! ] \varnothing, \, x]\! ]} \ee^{V(z)}} \le \frac1r \ee^{-V(x)}$ for $x\in \mathscr{L}_r$; hence
\begin{equation}
    P_\omega \Big( \bigcup_{i=1}^n\{ X_i \in \mathscr{L}_r\} \Big)
    \le
    \frac{E_\omega (L_n(\varnothing)+1)}{r} \sum_{x\in \mathscr{L}_r}\ee^{-V(x)} .
    \label{pf:l:main}
\end{equation}

We use the trivial inequality $L_n(\varnothing) \le n$, so $E_\omega (L_n(\varnothing)) \le n$. We now take $r=r(n) := \frac{n}{(\log n)^\gamma}$. With this choice of $r$, Lemma \ref{l:sum(xinLn)} tells us that if $\gamma<1$, then $(\log n)^\gamma \sum_{x\in \mathscr{L}_r}\ee^{-V(x)} \to 0$ in $\P^*$-probability. This yields Lemma \ref{l:ligne-arret}.\hfill$\Box$

\subsection{Covariance for edge local times of biased walks}
\label{subs:covariance}

$\phantom{aob}$In the proof of Proposition \ref{p:x=0} in Section \ref{s:proof-p:x=0}, we are going to estimate the covariance of local time of the biased walk $(X_i)$. It turns out to be more convenient to deal with covariance of {\it edge} local time instead of site local time. More precisely, for any $k\ge 1$ and any vertex $x\in \T\backslash\{\varnothing\}$, let us define the edge local time
\begin{equation}
    \overline{L}_k(x)
    :=
    \sum_{i=1}^k {\bf 1}_{\{ X_{i-1} = {\buildrel \leftarrow \over x} , \, X_i =x\} } \, ,
    \label{edge_local_time}
\end{equation}

\noindent which is the number of passages of the walk $(X_i)$, in the first $k$ steps, on the oriented edge from ${\buildrel \leftarrow \over x}$ to $x$. We are interested in the (edge) local time during an excursion away from $\varnothing$. 

Notation: $x\wedge y$ is the youngest common ancestor of $x$ and $y$ (or, equivalently, the unique vertex satisfying $[\![ \varnothing, \, x\wedge y]\!] = [\![ \varnothing, \, x ]\!] \cap [\![ \varnothing, \, y]\!]$).

\medskip

\begin{lemma}
\label{l:Cov}

 Let $T_\varnothing^+ := \inf\{ i\ge 1: \, X_i = \varnothing \}$ denote the first return to the root $\varnothing$ as in $(\ref{tau})$.
 
 {\rm (i)} We have, for $x \not= y \in \T$,
 \begin{equation}
     \mathrm{Cov}_\omega
     [\overline{L}_{T_\varnothing^+} (x), \, \overline{L}_{T_\varnothing^+} (y) ]
     \le 
     2\, \ee^{-[V(x)-V(x\wedge y)]-[V(y)-V(x\wedge y)]}\,
     E_\omega [\overline{L}_{T_\varnothing^+}(x\wedge y)^2] \, ,
     \label{covariance-negative}
 \end{equation}

 \noindent where $\mathrm{Cov}_\omega$ stands for covariance under the quenched probability $P_\omega$.
 
 {\rm (ii)} We have, for $x\in \T \backslash\{ \varnothing\}$,
 \begin{eqnarray}
     E_\omega [\overline{L}_{T_\varnothing^+}(x)] 
  &=& \omega(\varnothing, \, {\buildrel \leftarrow \over \varnothing}) \, \ee^{-V(x)} \, .
     \label{edge_lt_first_moment}
     \\
     E_\omega [\overline{L}_{T_\varnothing^+}(x)^2] 
  &=&\omega(\varnothing, \, {\buildrel \leftarrow \over \varnothing})\, \ee^{-V(x)} \Big( 2 \sum_{y\in \, ]\!] \varnothing, \, x]\!]} \ee^{V(y)-V(x)} -1 \Big) \, .
     \label{edge_lt_second_moment}
 \end{eqnarray}
  
\end{lemma}

\medskip

\noindent {\it Proof.} (i) We use the following elementary identity: for any pairs of random variables $\xi_1$ and $\xi_2$ defined on a probability space $(\Omega, \, \mathscr{F}, \, \p)$, having finite second moments, and any $\sigma$-field $\mathscr{G} \subset \mathscr{F}$, we have
\begin{equation}
    \mathrm{Cov}(\xi_1, \, \xi_2)
    =
    \e \Big[ \mathrm{Cov}_{\mathscr{G}} (\xi_1, \, \xi_2) \Big]
    +
    \mathrm{Cov} \Big[ \e(\xi_1\, | \, \mathscr{G}), \, \e(\xi_2\, | \, \mathscr{G}) \Big] \, ,
    \label{covar_conditionnelle}
\end{equation}

\noindent where $\mathrm{Cov}_{\mathscr{G}} (\xi_1, \, \xi_2) := \e(\xi_1 \, \xi_2 \, | \, \mathscr{G}) - \e(\xi_1\, | \, \mathscr{G}) \, \e(\xi_2\, | \, \mathscr{G})$ is the conditional covariance of $\xi_1$ and $\xi_2$ given $\mathscr{G}$.

We first treat the case that neither of $x$ and $y$ is an ancestor of the other. 

We write $u = u(x, \, y):= x\wedge y$ for brevity, and let $x_{|u|+1}$ and $y_{|u|+1}$ be the ancestor, at generation $|u|+1$, of $x$ and $y$ respectively. By definition of $x\wedge y$, the vertices $x_{|u|+1}$ and $y_{|u|+1}$ are distinct children of $u$. Conditionally on $\overline{L}_{T_\varnothing^+}(x_{|u|+1})$ and $\overline{L}_{T_\varnothing^+}(y_{|u|+1})$, the edge local times $\overline{L}_{T_\varnothing^+}(x)$ and $\overline{L}_{T_\varnothing^+}(y)$ are independent. We apply (\ref{covar_conditionnelle}) to $\xi_1 := \overline{L}_{T_\varnothing^+}(x)$, $\xi_2 := \overline{L}_{T_\varnothing^+}(y)$ and $\mathscr{G}:= \sigma(\overline{L}_{T_\varnothing^+}(x_{|u|+1}), \, \overline{L}_{T_\varnothing^+}(y_{|u|+1}))$, the $\sigma$-field generated by the edge local times $\overline{L}_{T_\varnothing^+}(x_{|u|+1})$ and $\overline{L}_{T_\varnothing^+}(y_{|u|+1})$. Since the conditional covariance vanishes, (\ref{covar_conditionnelle}) gives that
\begin{equation}
    \mathrm{Cov}_\omega [\overline{L}_{T_\varnothing^+}(x), \, \overline{L}_{T_\varnothing^+}(y)]
    =
    \mathrm{Cov}_\omega [ E_\omega (\overline{L}_{T_\varnothing^+}(x) \, | \, \mathscr{G}), \, E_\omega (\overline{L}_{T_\varnothing^+}(y) \, | \, \mathscr{G})] \, ,
    \label{covariance1}
\end{equation}

\noindent with $\mathscr{G}:= \sigma(\overline{L}_{T_\varnothing^+}(x_{|u|+1}), \, \overline{L}_{T_\varnothing^+}(y_{|u|+1}))$. Let us compute $E_\omega (\overline{L}_{T_\varnothing^+}(x) \, | \, \mathscr{G})$, which is nothing else but $E_\omega (\overline{L}_{T_\varnothing^+}(x) \, | \, \overline{L}_{T_\varnothing^+}(x_{|u|+1}))$. Write $|x| =: j > i:= |u|$. Then for any $k\in (i, \, j)\cap \z$, and given $\overline{L}_{T_\varnothing^+}(x_k) = \ell \ge 1$, $\overline{L}_{T_\varnothing^+}(x_{k+1})$ has the law of $\sum_{m=1}^\ell G_m$, where $G_m$, $m\ge 1$, are i.i.d.\ geometric random variables with parameter $p_k := \frac{\omega(x_k, \, x_{k-1})}{\omega(x_k, \, x_{k+1}) + \omega(x_k, \, x_{k-1})}$ (i.e., $G_m$ takes value $r$ with probability $(1-p_k)^rp_k$ for all non-negative integer $r$). Since $G_m$ has mean $\frac{1-p_k}{p_k}$, we have $E_\omega (\overline{L}_{T_\varnothing^+}(x_{k+1}) \, | \, \overline{L}_{T_\varnothing^+}(x_k)) = \overline{L}_{T_\varnothing^+}(x_k)\, \frac{1-p_k}{p_k} = \overline{L}_{T_\varnothing^+}(x_k) \, \ee^{-[V(x_{k+1})-V(x_k)]}$. As a consequence, we deduce from the Markov property of $k \to  \overline{L}_{T_\varnothing^+}(x_k) $ (under $P_\omega$) that 
\begin{eqnarray*}
    E_\omega (\overline{L}_{T_\varnothing^+}(x) \, | \, \overline{L}_{T_\varnothing^+}(x_{|u|+1}))
 &=&\overline{L}_{T_\varnothing^+}(x_{|u|+1}) \prod_{k=i+1}^{j-1} \ee^{-[V(x_{k+1})-V(x_k)]}
    \\
 &=&\overline{L}_{T_\varnothing^+}(x_{|u|+1}) \, \ee^{-[V(x)-V(x_{|u|+1})]} \, .
\end{eqnarray*}

\noindent Similarly, $E_\omega (\overline{L}_{T_\varnothing^+}(y) \, | \, \overline{L}_{T_\varnothing^+}(y_{|u|+1})) = \overline{L}_{T_\varnothing^+}(y_{|u|+1}) \, \ee^{-[V(y)-V(y_{|u|+1})]}$. Going back to (\ref{covariance1}), this leads to:
\begin{eqnarray}
 &&\mathrm{Cov}_\omega [\overline{L}_{T_\varnothing^+}(x), \, \overline{L}_{T_\varnothing^+}(y)]
    \nonumber
    \\
 &=&\ee^{-[V(x)-V(x_{|u|+1})] -[V(y)-V(y_{|u|+1})]} \, \mathrm{Cov}_\omega [\overline{L}_{T_\varnothing^+}(x_{|u|+1}), \, \overline{L}_{T_\varnothing^+}(y_{|u|+1}) ] \, .
    \label{covariance2}
\end{eqnarray}

\noindent To compute the covariance on the right-hand side, we write $(u^{(1)}, \, \cdots , \, u^{(N(u))})$ for the children of $u$ (among which are $x_{|u|+1}$ and $y_{|u|+1}$; so $N(u) \ge 2$), and observe that conditionally on $\overline{L}_{T_\varnothing^+}(u) = \ell \ge 1$, the law of the random vector $(\overline{L}_{T_\varnothing^+}(u^{(k)}), \, 1\le k\le N(u))$ under $P_\omega$ is multinomial with parameter $(\sum_{k=1}^\ell \mathfrak{G}_k, \, (p^{(k)}(u) := \frac{\omega(u, \, u^{(k)})}{1- \omega(u, \, {\buildrel \leftarrow \over u})}, \, 1\le k\le N(u)))$, where $\mathfrak{G}_k$, $k\ge 1$, are i.i.d.\ random variables having the geometric distribution of parameter $\omega(u, \, {\buildrel \leftarrow \over u})$.\footnote{A random vector $(\xi_1, \, \cdots, \, \xi_N)$ defined on $(\Omega, \, \mathscr{F}, \, \p)$ has the multinomial distribution with parameter $(m, \, (p^{(1)}, \, \cdots, \, p^{(N)}))$ if $\p(\xi_1 = m_1, \, \cdots, \, \xi_N=m_N) = \frac{m!}{m_1! \, \cdots \, m_N!} \prod_{k=1}^N (p^{(k)})^{m_k}$ for all non-negative integers $m_k$, $1\le k\le N$, satisfying $m_1+ \cdots m_N=m$; in particular, $\e (s_1^{\xi_1} \cdots s_N^{\xi_N}) = (\sum_{k=1}^N p^{(k)} s_k)^m$, for all $s_k \ge 0$, $1\le k\le N$.} Accordingly, for all $\ell \ge 1$,
\begin{eqnarray*}
    E_\omega \Big[ \prod_{k=1}^{N(u)} (s_k)^{\overline{L}_{T_\varnothing^+}(u^{(k)})} \, \Big| \, \overline{L}_{T_\varnothing^+}(u) = \ell \Big]
 &=& E_\omega \Big[ \Big( \sum_{k=1}^{N(u)} \frac{s_k \, \omega(u, \, u^{(k)})}{1- \omega(u, \, {\buildrel \leftarrow \over u})} \Big)^{\sum_{k=1}^\ell \mathfrak{G}_k} \, \Big]
    \\
 &=& \Big\{ E_\omega \Big[ \Big( \sum_{k=1}^{N(u)} \frac{s_k \, \omega(u, \, u^{(k)})}{1- \omega(u, \, {\buildrel \leftarrow \over u})} \Big)^{\mathfrak{G}_1}\, \Big] \Big\}^\ell .
\end{eqnarray*}
 
\noindent Since $E_\omega (s^{\mathfrak{G}_1}) = \frac{ \omega(u, \, {\buildrel \leftarrow \over u})}{1- s (1- \omega(u, \, {\buildrel \leftarrow \over u}))}$, this yields
\begin{equation}
    E_\omega \Big[ \prod_{k=1}^{N(u)} (s_k)^{\overline{L}_{T_\varnothing^+}(u^{(k)})} \, \Big| \, \overline{L}_{T_\varnothing^+}(u) \Big]
    =
    \Big\{   \frac{ \omega(u, \, {\buildrel \leftarrow \over u}) }{   1- \sum_{k=1}^{N(u)} s_k \, \omega(u, \, u^{(k)})} \Big\}^{\overline{L}_{T_\varnothing^+}(u)} \, .
    \label{fct_generatrice_edge_lt}
\end{equation}

\noindent [We proved it assuming that $\overline{L}_{T_\varnothing^+}(u) \ge 1$, but it is trivially true on the set $\{ \overline{L}_{T_\varnothing^+}(u) =0\}$.] In particular, for $1\le k\not= m \le N(u)$,
\begin{eqnarray*}
    E_\omega [\overline{L}_{T_\varnothing^+}(u^{(k)}) \, | \, \overline{L}_{T_\varnothing^+}(u)]
 &=&\ee^{-[V(u^{(k)})-V(u)]} \, \overline{L}_{T_\varnothing^+}(u)\, ,
    \\
    E_\omega [\overline{L}_{T_\varnothing^+}(u^{(k)})\, \overline{L}_{T_\varnothing^+}(u^{(m)})  \, | \, \overline{L}_{T_\varnothing^+}(u)]
 &=& \ee^{-[V(u^{(k)})-V(u)]-[V(u^{(m)})-V(u)]} \, 
    \overline{L}_{T_\varnothing^+}(u) \,
    (\overline{L}_{T_\varnothing^+}(u)+1) \, .
\end{eqnarray*}

\noindent Applying again (\ref{covar_conditionnelle}), this time to $\xi_1 := \overline{L}_{T_\varnothing^+}(u^{(k)})$, $\xi_2 := \overline{L}_{T_\varnothing^+}(u^{(m)})$ and $\mathscr{G}:= \sigma( \overline{L}_{T_\varnothing^+}(u))$, we obtain ($\mathrm{Var}_\omega$ denoting variance under $P_\omega$): 
\begin{eqnarray*}
    \mathrm{Cov}_\omega [\overline{L}_{T_\varnothing^+}(u^{(k)}), \, \overline{L}_{T_\varnothing^+}(u^{(m)}) ] 
 &=& \ee^{-[V(u^{(k)})-V(u)]-[V(u^{(m)})-V(u)]} \, E_\omega [\overline{L}_{T_\varnothing^+}(u)]
    \\
 &&+\ee^{-[V(u^{(k)})-V(u)]-[V(u^{(m)})-V(u)]} \, \mathrm{Var}_\omega [\overline{L}_{T_\varnothing^+}(u)]
    \\
 &\le&2\, \ee^{-[V(u^{(k)})-V(u)]-[V(u^{(m)})-V(u)]} \, E_\omega [\overline{L}_{T_\varnothing^+}(u)^2] \, ,
\end{eqnarray*}

\noindent the last inequality following from the fact that $\overline{L}_{T_\varnothing^+}(u) \le \overline{L}_{T_\varnothing^+}(u)^2$ (recalling that $\overline{L}_{T_\varnothing^+}(u)$ is integer-valued). We take $k$ and $m$ be such that $u^{(k)} = x_{|u|+1}$ and $u^{(m)} = y_{|u|+1}$. In view of (\ref{covariance2}), this yields the desired inequality \eqref{covariance-negative} in Lemma \ref{l:Cov}.


It remains to deal with the special case that either $x$ is an ancestor of $y$, or $y$ is an ancestor of $x$. 

This, however, is easy. Without loss of generality, let us assume that $y$ is an ancestor of $x$, in which case we have seen that $E_\omega (\overline{L}_{T_\varnothing^+}(x) \, | \, \overline{L}_{T_\varnothing^+}(y)) = \ee^{-[V(x)-V(y)]}\, \overline{L}_{T_\varnothing^+}(y)$. So applying (\ref{covar_conditionnelle}) to $\xi_1 := \overline{L}_{T_\varnothing^+}(x)$, $\xi_2 := \overline{L}_{T_\varnothing^+}(y)$ and $\mathscr{G} := \sigma( \overline{L}_{T_\varnothing^+}(y) )$ gives
$$
\mathrm{Cov}_\omega [\overline{L}_{T_\varnothing^+}(x), \, \overline{L}_{T_\varnothing^+}(y) ] 
=
0
+
\ee^{-[V(x)-V(y)]}\, \mathrm{Var}_\omega [\overline{L}_{T_\varnothing^+}(y)] \, ,
$$

\noindent yielding \eqref{covariance-negative}. 

(ii) We already noted that $E_\omega [\overline{L}_{T_\varnothing^+}(x)] = \ee^{-[V(x)-V(x_1)]} E_\omega [ \overline{L}_{T_\varnothing^+}(x_1)]$, where $x_1$ denotes, as before, the ancestor of $x$ in the first generation. Since $E_\omega [ \overline{L}_{T_\varnothing^+}(x_1)] = \omega(\varnothing, \, x_1)$, and by definition, $\omega(\varnothing, \, x_1) = \omega(\varnothing, \, {\buildrel \leftarrow \over \varnothing}) \, \ee^{-V(x_1)}$, this yields $E_\omega [\overline{L}_{T_\varnothing^+}(x)] = \omega(\varnothing, \, {\buildrel \leftarrow \over \varnothing}) \, \ee^{-V(x)}$, as stated in \eqref{edge_lt_first_moment}.
 
It remains to compute $E_\omega [\overline{L}_{T_\varnothing^+}(x)^2]$. {F}rom \eqref{fct_generatrice_edge_lt}, we get that
$$
E_\omega [\overline{L}_{T_\varnothing^+}(u^{(k)}) [\overline{L}_{T_\varnothing^+}(u^{(k)})-1] \, | \, \overline{L}_{T_\varnothing^+}(u)]
=
\ee^{-2[V(u^{(k)})-V(u)]} \, \overline{L}_{T_\varnothing^+}(u) [\overline{L}_{T_\varnothing^+}(u)+1] \, .
$$

\noindent Taking expectation on both sides, and replacing the pair $(u^{(k)}, \, u)$ by $(x, \, {\buildrel \leftarrow \over x})$, we obtain:
$$
E_\omega [\overline{L}_{T_\varnothing^+}(x)^2]
=
\ee^{-2[V(x)-V({\buildrel \leftarrow \over x})]} \, E_\omega [\overline{L}_{T_\varnothing^+}({\buildrel \leftarrow \over x})^2]
+
(\ee^{-2[V(x)-V({\buildrel \leftarrow \over x})]} + \ee^{-[V(x)-V({\buildrel \leftarrow \over x})]})\, E_\omega [\overline{L}_{T_\varnothing^+}({\buildrel \leftarrow \over x})]
$$

\noindent By the already proved \eqref{edge_lt_first_moment}, $E_\omega [\overline{L}_{T_\varnothing^+}({\buildrel \leftarrow \over x})] = \omega(\varnothing, \, {\buildrel \leftarrow \over \varnothing}) \, \ee^{-V({\buildrel \leftarrow \over x})}$. Solving this difference equation (with initial condition $E_\omega [ \overline{L}_{T_\varnothing^+}(x_1)^2] = \omega(\varnothing, \, x_1) = \omega(\varnothing, \, {\buildrel \leftarrow \over \varnothing}) \, \ee^{-V(x_1)}$) yields \eqref{edge_lt_second_moment}. This completes the proof of the lemma.\hfill$\Box$

\section{Biased walks: proof of Proposition \ref{p:x=0}}
\label{s:proof-p:x=0}

$\phantom{aob}$Let $P_\omega^{(r)}$ denote the quenched law of the biased walk {\it with a reflecting barrier at} $\mathscr{L}_r$. 
Under $P_\omega^{(r)}$, the biased walk $(X_i)$ is positive recurrent taking values in $\{ x\in \T: \, x\le \mathscr{L}_r\} \cup \{ {\buildrel \leftarrow \over \varnothing} \}$, with invariant probability $\pi_r$ as in (\ref{pi_n}). In particular, if $T_\varnothing^+$ denotes, as in (\ref{tau}), the first return time to $\varnothing$, and $L_{T_\varnothing^+}$ (site) local time as in \eqref{local_time}, 
\begin{eqnarray}
    E_\omega^{(r)} (T_\varnothing^+)
 &=& \frac{1}{\pi_r(\varnothing)}\, ,
    \label{E(T)}
    \\
    E_\omega^{(r)} [L_{T_\varnothing^+} (y) ]
 &=&\frac{\pi_r(y)}{\pi_r(\varnothing)}\, ,
    \qquad
    y\in \{ x\in \T: \, x\le \mathscr{L}_r\} \cup \{ {\buildrel \leftarrow \over \varnothing} \}  \, .
    \label{E(L)}
\end{eqnarray}

We now proceed to study $L_n (\varnothing)$ under $P_\omega^{(r)}$. Let $\ell \ge 1$, and let $T_\varnothing^{(\ell)}$ denote the $\ell$-th return time to $\varnothing$ (so $T_\varnothing^{(1)}$ is $T_\varnothing^+$, under $P_\omega^{(r)}$). Under $P_\omega^{(r)}$, $T_\varnothing^{(\ell)}$ is the sum of $\ell$ independent copies of $T_\varnothing^+$. In particular, $E_\omega^{(r)} (T_\varnothing^{(\ell)}) = \ell \times E_\omega^{(r)} (T_\varnothing^+) = \frac{\ell}{\pi_r(\varnothing)}$.

By the simple relation $\{ L_n (\varnothing) \le \ell \} = \{ T_\varnothing^{(\ell)} \ge n \}$, we have
$$
    P_\omega^{(r)} \{ L_n (\varnothing) \le \ell \}
    =
    P_\omega^{(r)} \Big\{ T_\varnothing^{(\ell)} - \frac{\ell}{\pi_r(\varnothing)} \ge n- \frac{\ell}{\pi_r(\varnothing)} \Big\} \, ,
$$

\noindent which, by Chebyshev's inequality, is bounded by $\frac{\ell}{(n- \frac{\ell}{\pi_r(\varnothing)})^2} \, \mathrm{Var}_\omega^{(r)} (T_\varnothing^+)$ if $n>\frac{\ell}{\pi_r(\varnothing)}$ ($\mathrm{Var}_\omega^{(r)}$ denoting the variance under the probability $P_\omega^{(r)}$). However, it has not been clear to us whether $\mathrm{Var}_\omega^{(r)} (T_\varnothing^+)$ is sufficiently small. This is why some care is in order when applying the method of second moment. We are not going to estimate the variance (under $P_\omega^{(r)}$) of $T_\varnothing^+$; instead, we are going to decompose $T_\varnothing^+$ into three distinct parts, in such a way that the variance of a part is sufficiently small for our needs and that the expectation of the other parts is also sufficiently small. 

Recall from (\ref{gamma}) that $\mathscr{L}_r := \{ x: \, \sum_{z\in \, ]\! ] \varnothing, \, x]\! ]} \ee^{V(z)-V(x)} > r , \; \sum_{z\in \, ]\! ]\varnothing, \, y]\! ]} \ee^{V(z)-V(y)} \le r, \; \forall y\in \, ]\!] \varnothing, \, x[\![ \, \}$. 
The reason for which we have not been able to make $\mathrm{Var}_\omega^{(r)} (T_\varnothing^+)$ small is that $r$ is too large.
Our solution is to consider {\bf two} scales: $\mathscr{L}_r$ and $\mathscr{L}_s$ with $s := \frac{r}{(\log r)^\theta} \le r$, where $\theta\ge 0$.

The promised decomposition for $T_\varnothing^+$
is as follows, the constant $\delta_1$ being defined in (\ref{integrabilite_sous_Q}):
\begin{eqnarray}
    T_\varnothing^{\,\mathrm{(a)}}
 &:=& \sum_{y\in \T: \, y < \mathscr{L}_s} L_{T_\varnothing^+} (y) \, {\bf 1}_{\{ \min_{u\in \, [\! [ \varnothing, \, y]\! ]} \omega(u, \, {\buildrel \leftarrow \over u}) \ge (\log r)^{-6/\delta_1} \} } \, ,
    \label{T(close)}
    \\
    T_\varnothing^{\,\mathrm{(b)}}
 &:=& \sum_{y\in \T: \, y < \mathscr{L}_s} L_{T_\varnothing^+} (y) \, {\bf 1}_{\{ \min_{u\in \, [\! [ \varnothing, \, y]\! ]} \omega(u, \, {\buildrel \leftarrow \over u}) < (\log r)^{-6/\delta_1} \} }  \, ,
    \label{T(close,bad)}
    \\
    T_\varnothing^{\,\mathrm{(c)}}
 &:=& \sum_{y\in \T: \, \mathscr{L}_s\le y \le \mathscr{L}_r} L_{T_\varnothing^+} (y) \, .
    \label{T(far)}
\end{eqnarray}

\noindent Then
\begin{equation}
    T_\varnothing^+ -1 
    \le 
    T_\varnothing^{\,\mathrm{(a)}} + T_\varnothing^{\,\mathrm{(b)}} + T_\varnothing^{\,\mathrm{(c)}} 
    \le 
    T_\varnothing^+ \, .
    \label{T=Ta+Tb+Tc}
\end{equation}

\noindent [The quantities $T_\varnothing^+$ and $T_\varnothing^{\,\mathrm{(a)}} + T_\varnothing^{\,\mathrm{(b)}} + T_\varnothing^{\,\mathrm{(c)}}$ can differ by $1$ in case $X_1 = {\buildrel \leftarrow \over \varnothing}$.]

The next pair of lemmas summarize basic properties of $T_\varnothing^{\,\mathrm{(a)}}$, $T_\varnothing^{\,\mathrm{(b)}}$ and $T_\varnothing^{\,\mathrm{(c)}}$ that are needed in this paper: loosely speaking, we control in a satisfying way the first two moments of $T_\varnothing^{\,\mathrm{(a)}}$, and although we have no control on the variances of $T_\varnothing^{\,\mathrm{(b)}}$ and $T_\varnothing^{\,\mathrm{(c)}}$, we show that they both have negligible expectations compared to the expectation of $T_\varnothing^+$.
 
\medskip

\begin{lemma}
\label{l:E(T)}

 Let $\theta\ge 0$ and let $s := \frac{r}{(\log r)^\theta}$. When $r\to \infty$,
 \begin{eqnarray}
     \frac{E_\omega^{(r)} (T_\varnothing^{\,\mathrm{(a)}})}{E_\omega^{(r)} (T_\varnothing^+)}
  &\to& 1,
     \qquad \hbox{in $\P^*$-probability},
     \label{E(T(close))}
     \\
     \frac{E_\omega^{(r)} (T_\varnothing^{\,\mathrm{(b)}})}{E_\omega^{(r)} (T_\varnothing^+)}
  &\to& 0,
     \qquad \hbox{in $\P^*$-probability},
     \label{E(T(close,bad))}
     \\
     \frac{E_\omega^{(r)} (T_\varnothing^{\,\mathrm{(c)}})}{E_\omega^{(r)} (T_\varnothing^+)}
  &\to& 0,
     \qquad \hbox{in $\P^*$-probability}.
     \label{E(T(far))}
 \end{eqnarray}

 \noindent In particular,
 \begin{equation}
     \frac{E_\omega^{(r)} (T_\varnothing^{\,\mathrm{(b)}}) + E_\omega^{(r)} (T_\varnothing^{\,\mathrm{(c)}})}{\log r}
     \to
     0, 
     \qquad \hbox{in $\P^*$-probability}.
     \label{E(Tb)+E(Tc)=o(logr)}
 \end{equation}

\end{lemma}

\medskip

\begin{lemma}
\label{l:V(T)}

 Let $\theta\ge 0$ and let $s := \frac{r}{(\log r)^\theta}$. There exists a constant $c_7>0$ such that for all $r\ge 2$,
 \begin{equation}
     \E \Big[ \, \mathrm{Var}_\omega^{(r)} (T_\varnothing^{\,\mathrm{(a)}})\, \Big]
     \le
     c_7 \, s \, (\log r)^{\frac{18}{\delta_1}+6} \, ,
     \label{Var(T(close))}
 \end{equation}

 \noindent where $\delta_1>0$ is the constant in $(\ref{integrabilite_sous_Q})$. 

\end{lemma}

\medskip

By admitting Lemmas \ref{l:E(T)} and \ref{l:V(T)} for the time being, we are able to prove Proposition \ref{p:x=0}.

\bigskip

\noindent {\it Proof of Proposition \ref{p:x=0}.} Let $\theta\ge 0$ and let $s := \frac{r}{(\log r)^\theta}$. Let
\begin{eqnarray*}
    T_\varnothing^{(\ell), \,\mathrm{(a)}}
 &:=&\sum_{y\in \T: \, y < \mathscr{L}_s} L_{T_\varnothing^{(\ell)}} (y) \, {\bf 1}_{\{ \min_{u\in \, [\! [ \varnothing, \, y]\! ]} \omega(u, \, {\buildrel \leftarrow \over u}) \ge (\log r)^{-6/\delta_1} \} } \, ,
    \\
    T_\varnothing^{(\ell), \,\mathrm{(b)}}
 &:=&\sum_{y\in \T: \, y < \mathscr{L}_s} L_{T_\varnothing^{(\ell)}} (y) \, {\bf 1}_{\{ \min_{u\in \, [\! [ \varnothing, \, y]\! ]} \omega(u, \, {\buildrel \leftarrow \over u}) < (\log r)^{-6/\delta_1} \} } \, ,
    \\
    T_\varnothing^{(\ell), \,\mathrm{(c)}}
 &:=&\sum_{y\in \T: \, \mathscr{L}_s \le y \le \mathscr{L}_r} L_{T_\varnothing^{(\ell)}} (y) \, .
\end{eqnarray*}

\noindent Then $T_\varnothing^{(\ell)} - \ell \le T_\varnothing^{(\ell), \,\mathrm{(a)}} + T_\varnothing^{(\ell), \,\mathrm{(b)}} + T_\varnothing^{(\ell), \,\mathrm{(c)}} \le T_\varnothing^{(\ell)}$.

For any $n_1\ge 1$ and $n_2 \ge 1$ with $n_1+n_2\le n-\ell$,
\begin{eqnarray*}
    P_\omega^{(r)} \{ L_n (\varnothing) \le \ell \}
 &=& P_\omega^{(r)} \{ T_\varnothing^{(\ell)} \ge n\}
    \\
 &\le&P_\omega^{(r)} \{ T_\varnothing^{(\ell), \,\mathrm{(a)}} + T_\varnothing^{(\ell), \,\mathrm{(b)}} + T_\varnothing^{(\ell), \,\mathrm{(c)}} \ge n-\ell\}
    \\
 &\le&P_\omega^{(r)} \{ T_\varnothing^{(\ell), \,\mathrm{(a)}} \ge n_1\}
    +
    P_\omega^{(r)} \{ T_\varnothing^{(\ell), \,\mathrm{(b)}} + T_\varnothing^{(\ell), \,\mathrm{(c)}} \ge n_2\} \, .
\end{eqnarray*}

\noindent Observe that $E_\omega^{(r)} [T_\varnothing^{(\ell), \,\mathrm{(b)}} + T_\varnothing^{(\ell), \,\mathrm{(c)}}] = \ell \, [ E_\omega^{(r)} (T_\varnothing^{\,\mathrm{(b)}}) + E_\omega^{(r)} (T_\varnothing^{\,\mathrm{(c)}})]$, so by Markov's inequality,
$$
P_\omega^{(r)} \{ T_\varnothing^{(\ell), \,\mathrm{(b)}} + T_\varnothing^{(\ell), \,\mathrm{(c)}} \ge n_2\}
\le
\frac{\ell}{n_2} \, [ E_\omega^{(r)} (T_\varnothing^{\,\mathrm{(b)}}) + E_\omega^{(r)} (T_\varnothing^{\,\mathrm{(c)}})] \, .
$$

\noindent For $P_\omega^{(r)} \{ T_\varnothing^{(\ell), \,\mathrm{(a)}} \ge n_1\}$, we note that $E_\omega^{(r)} (T_\varnothing^{(\ell), \,\mathrm{(a)}}) = \ell \times E_\omega^{(r)} (T_\varnothing^{\,\mathrm{(a)}}) \le \ell \times E_\omega^{(r)} (T_\varnothing^+) = \frac{\ell}{\pi_r(\varnothing)}$, and that $\mathrm{Var}_\omega^{(r)} (T_\varnothing^{(\ell), \,\mathrm{(a)}}) = \ell \, \mathrm{Var}_\omega^{(r)} (T_\varnothing^{\,\mathrm{(a)}})$. If $n_1 - \frac{\ell}{\pi_r(\varnothing)}>0$, then by Chebyshev's inequality,
\begin{eqnarray*}
    P_\omega^{(r)} \{ T_\varnothing^{(\ell), \,\mathrm{(a)}} \ge n_1\}
 &\le&P_\omega^{(r)} \{ T_\varnothing^{(\ell), \,\mathrm{(a)}} - E_\omega^{(r)} (T_\varnothing^{(\ell), \,\mathrm{(a)}}) \ge n_1 - \frac{\ell}{\pi_r(\varnothing)} \}
    \\
 &\le&\frac{\ell \, \mathrm{Var}_\omega^{(r)} (T_\varnothing^{\,\mathrm{(a)}})}{[n_1 - \frac{\ell}{\pi_r(\varnothing)}]^2} \, .
\end{eqnarray*}

\noindent Let us now fix the choice for $\ell$, $n_1$ and $n_2$. Let $0<\varepsilon<1$. We take $n_1 := \lceil (1+\varepsilon) \frac{\ell}{\pi_r(\varnothing)}\rceil$ and $n_2 := \lfloor \varepsilon \frac{\ell}{\pi_r(\varnothing)} \rfloor -\ell -1$ so that $n_1 + n_2 \le (1+2\varepsilon) \frac{\ell}{\pi_r(\varnothing)} -\ell$, which is indeed bounded by $n-\ell$ if we take $\ell := \lfloor \frac{1}{1+2\varepsilon}\, n\, \pi_r(\varnothing)\rfloor$. With the choice made for $(\ell, \, n_1, \, n_2)$, we have 
$$
P_\omega^{(r)} \{ L_n (\varnothing) \le \ell \}
\le
\frac{\ell \, \mathrm{Var}_\omega^{(r)} (T_\varnothing^{\,\mathrm{(a)}})}{[n_1 - \frac{\ell}{\pi_r(\varnothing)}]^2} 
+
\frac{\ell}{n_2} \, [ E_\omega^{(r)} (T_\varnothing^{\,\mathrm{(b)}}) + E_\omega^{(r)} (T_\varnothing^{\,\mathrm{(c)}})] \, .
$$

\noindent Recall 
from Theorem \ref{t:Yn} that $(\log r) \, \pi_r(\varnothing) \to \frac{\sigma^2}{4D_\infty} \, \ee^{-U(\varnothing)}$ in $\P^*$-probability (for $r\to \infty$). We choose $r=n$ so that we are entitled to apply Lemma \ref{l:ligne-arret}. With the definition of $s:= \frac{r}{(\log r)^\theta}$, we apply Lemma \ref{l:V(T)} (choosing $\theta > \frac{18}{\delta_1}+ 5$) and Lemma \ref{l:E(T)} (part \eqref{E(Tb)+E(Tc)=o(logr)}), to see that $P_\omega^{(r)} \{ L_n (\varnothing) \le \ell \} \to 0$ in $\P^*$-probability. By Lemma \ref{l:ligne-arret}, $\p ( \cup_{i=1}^n\{ X_i \in \mathscr{L}_r \} ) \to 0$, so this is equivalent to saying that $P_\omega \{ L_n (\varnothing) \le \ell \} \to 0$ in $\P^*$-probability, with the choice of $\ell := \lfloor \frac{1}{1+2\varepsilon}\, n\, \pi_r(\varnothing)\rfloor$. Again, since $(\log r) \, \pi_r(\varnothing) \to \frac{\sigma^2}{4D_\infty} \, \ee^{-U(\varnothing)}$ in $\P^*$-probability (Theorem \ref{t:Yn}), this yields the lower bound in \eqref{localtimeproba}.

The proof of the upper bound is similar, with the same choice $r:=n$, and is slightly easier because we do not need to care about $  T_\varnothing^{(\ell), \,\mathrm{(b)}}$ and $T_\varnothing^{(\ell), \,\mathrm{(c)}}$ any more. Indeed, for any $\ell \ge 1$, 
\begin{eqnarray*}
    P_\omega^{(r)} \{ L_n (\varnothing) \ge \ell \}
 &=&P_\omega^{(r)} \{ T_\varnothing^{(\ell)} \le n \}
    \\
 &\le&P_\omega^{(r)} \{ T_\varnothing^{(\ell), \,\mathrm{(a)}} \le n \}
    \\
 &\le&\frac{\mathrm{Var}_\omega^{(r)} (T_\varnothing^{(\ell), \,\mathrm{(a)}})}{[E_\omega^{(r)} (T_\varnothing^{(\ell), \,\mathrm{(a)}}) - n]^2} \, ,
\end{eqnarray*}

\noindent as long as $E_\omega^{(r)} (T_\varnothing^{(\ell), \,\mathrm{(a)}}) > n$. Again, $E_\omega^{(r)} (T_\varnothing^{(\ell), \,\mathrm{(a)}}) = \ell \, E_\omega^{(r)} (T_\varnothing^{\,\mathrm{(a)}})$, and $\mathrm{Var}_\omega^{(r)} (T_\varnothing^{(\ell), \,\mathrm{(a)}}) = \ell \, \mathrm{Var}_\omega^{(r)} (T_\varnothing^{\mathrm{(a)}})$. This time, with $\varepsilon>0$, our choice is $\ell := \lfloor (1+\varepsilon) \, n \, \pi_r(\varnothing)\rfloor$. For this new choice of $\ell$, part \eqref{E(T(close))} of Lemma \ref{l:E(T)} ensures that $\P^*\{ E_\omega^{(r)} (T_\varnothing^{(\ell), \,\mathrm{(a)}}) > (1+ \frac{\varepsilon}{2}) n\} \to 1$ for $n\to \infty$. So by Lemma \ref{l:V(T)}, if $s:= \frac{r}{(\log r)^\theta}$ with $\theta > \frac{18}{\delta_1}+ 5$, then $P_\omega^{(r)} \{ L_n (\varnothing) \ge \ell \} \to 0$ in $\P^*$-probability, which, in view of Lemma \ref{l:ligne-arret}, is equivalent to saying that $P_\omega \{ L_n (\varnothing) \ge \ell \} \to 0$ in $\P^*$-probability. This yields \eqref{localtimeproba}.

It remains to check \eqref{localtimemoment}. In view of \eqref{localtimeproba}, it suffices to show the following:
\begin{equation}
    \Big( \frac{\log n}{n}\Big)^2 \, E_\omega [ (L_n(\varnothing))^2 ]
    \hbox{ \rm is tight under } \P^* \, .
    \label{E(L)<n/log(n)}
\end{equation}

Clearly, $E_\omega [(L_n(\varnothing))^2] \le E_\omega^{(r)} [ (L_n(\varnothing))^2]$, for any $r>1$. Observe that 
$$
E_\omega^{(r)} [(L_n(\varnothing))^2]
\le
2 \sum_{j=1}^\infty j \, P_\omega^{(r)} \{ L_n(\varnothing) \ge j\}
=
2 \sum_{j=1}^\infty j\, P_\omega^{(r)} \{ T^{(j)}_\varnothing \le n\} \, .
$$

\noindent By Chebyshev's inequality, $P_\omega^{(r)} \{ T^{(j)}_\varnothing \le n\} \le \ee\times E_\omega^{(r)} (\ee^{- T^{(j)}_\varnothing/n})$, which, by the strong Markov property, is $\ee\times [E_\omega^{(r)} (\ee^{- T^+_\varnothing/n})]^j$. As such,
$$
E_\omega [(L_n(\varnothing))^2]
\le
2 \ee \, \sum_{j=1}^\infty\,  j\, \Big[ E_\omega^{(r)} (\ee^{- T^+_\varnothing/n} ) \Big]^j 
\le
\frac{2 \ee}{[ 1-  E_\omega^{(r)} (\ee^{- T^+_\varnothing/n}) ]^2},
$$
  
\noindent where, in the last inequality, we used the elementary fact that $\sum_{j=1}^\infty j x^j = \frac{x}{(1-x)^2} \le \frac{1}{(1-x)^2}$ for any $x \in [0, \, 1]$. 

Note that for any nonnegative random variable $\xi$ with $\e(\xi^2)<\infty$, we have $\e (1- \ee^{-\xi}) \ge \e (\xi - \frac{\xi^2}{2}) = \e (\xi)- \frac12 [\e(\xi)]^2 - \frac12 \mathrm{Var} (\xi)$. Therefore, 
\begin{eqnarray*}
    1 - E_\omega^{(r)} (\ee^{- T^+_\varnothing/n})
 &\ge& 1 - E_\omega^{(r)} (\ee^{- T^{(a)}_\varnothing/n})
    \\
 &\ge& \frac{E_\omega^{(r)} (T^{(a)}_\varnothing)}{n}
    - 
    \frac{[E_\omega^{(r)} (T^{(a)}_\varnothing )]^2}{2n^2} 
    -
    \frac{\mathrm{Var}_\omega^{(r)}(T^{(a)}_\varnothing)}{2n^2}\, .
\end{eqnarray*}
    
\noindent We choose again $r:= n$. By Lemma \ref{l:E(T)} and Theorem \ref{t:Yn}, $E_\omega^{(r)} (T^{(a)}_\varnothing) = (\frac{4}{\sigma^2} \, D_\infty + o_{\P^*}(1)) \log r$, where $o_{\P^*}(1)$ denotes, as before, a term converging to 0 in $\P^*$-probability and its value may vary from line to line. By Lemma \ref{l:V(T)}, if we choose $s:= \frac{r}{(\log r)^\theta}$ with $\theta > \frac{18}{\delta_1} +5$, then $\frac{\mathrm{Var}_\omega^{(r)} (T^{(a)}_\varnothing)}{n\log n} \to 0$ in $\P^*$-probability. This yields \eqref{E(L)<n/log(n)}, and thus \eqref{localtimemoment}. Proposition \ref{p:x=0} is proved.\hfill$\Box$



\bigskip

The rest of the section is devoted to the proof of Lemmas \ref{l:E(T)} and \ref{l:V(T)}.

\bigskip

\noindent {\it Proof of Lemma \ref{l:E(T)}.} Clearly, \eqref{E(Tb)+E(Tc)=o(logr)} follows from \eqref{E(T(close,bad))} and \eqref{E(T(far))} (combined with Theorem \ref{t:Yn}). On the other hand, $T_\varnothing^{\,\mathrm{(a)}} + T_\varnothing^{\,\mathrm{(b)}} + T_\varnothing^{\,\mathrm{(c)}}$ and $T_\varnothing^+$ differ by at most $1$ (see \eqref{T=Ta+Tb+Tc}), so (\ref{E(T(close,bad))}) and (\ref{E(T(far))}) together imply (\ref{E(T(close))}). As a consequence, we only need to prove \eqref{E(T(close,bad))} and \eqref{E(T(far))}.

Let us start with the {\bf proof of (\ref{E(T(far))})}. By definition of $T_\varnothing^{\,\mathrm{(c)}}$ (see \eqref{T(far)}), $E_\omega^{(r)} (T_\varnothing^{\,\mathrm{(c)}}) = \sum_{\mathscr{L}_s \le y\le \mathscr{L}_r} E_\omega^{(r)}[L_{T_\varnothing^+} (y)]$. We have seen in (\ref{E(L)}) that $E_\omega^{(r)} [L_{T_\varnothing^+} (y) ] = \frac{\pi_r(y)}{\pi_r(\varnothing)}$ for $y\le \mathscr{L}_r$. Since $E_\omega^{(r)} (T_\varnothing^+) = \frac{1}{\pi_r(\varnothing)}$ (see (\ref{E(T)})), we have 
$$
\frac{E_\omega^{(r)} (T_\varnothing^{\,\mathrm{(c)}})}{E_\omega^{(r)} (T_\varnothing^+)}
=
\sum_{\mathscr{L}_s \le y \le \mathscr{L}_r} \pi_r(y)
=
\sum_{\mathscr{L}_s \le y<\mathscr{L}_r} \frac{\ee^{-U(y)}}{Z_r} 
+
\sum_{y\in \mathscr{L}_r} \frac{\ee^{-V(y)}}{Z_r} 
 \, .
$$

\noindent Noting $\ee^{-U(y)} = \ee^{-V(y)} + \sum_{z\in \T: \, {\buildrel \leftarrow \over z} =y} \ee^{-V(z)}$, we arrive at:
$$
\frac{E_\omega^{(r)} (T_\varnothing^{\,\mathrm{(c)}})}{E_\omega^{(r)} (T_\varnothing^+)}
=
\sum_{\mathscr{L}_s \le y<\mathscr{L}_r} \frac{\ee^{-V(y)}}{Z_r} 
+
\sum_{\mathscr{L}_s <z\le\mathscr{L}_r} \frac{\ee^{-V(z)}}{Z_r}
+
\sum_{y\in \mathscr{L}_r} \frac{\ee^{-V(y)}}{Z_r} 
\le
\frac{2}{Z_r} \sum_{\mathscr{L}_s\le y\le \mathscr{L}_r} \ee^{-V(y)} \, .
$$

\noindent By Theorem \ref{t:Yn}, $\frac{1}{\log r}\sum_{y\le \mathscr{L}_r} \ee^{-V(y)} \to \frac{2}{\sigma^2} \, D_\infty$ in $\P^*$-probability and $\frac{1}{\log r}\sum_{y< \mathscr{L}_s} \ee^{-V(y)} \to \frac{2}{\sigma^2} \, D_\infty$ in $\P^*$-probability (noting that $\frac{1}{\log r} \sum_{y\in \mathscr{L}_s} \ee^{-V(y)} \to 0$ in $\P^*$-probability according to Lemma \ref{l:sum(xinLn)}). Hence $\frac{1}{\log r} \sum_{\mathscr{L}_s \le y\le \mathscr{L}_r} \ee^{-V(y)} \to \frac{2}{\sigma^2} \, D_\infty - \frac{2}{\sigma^2} \, D_\infty =0$ in $\P^*$-probability. Since $\frac{Z_r}{\log r} \to \frac{4}{\sigma^2} \, D_\infty >0$ in $\P^*$-probability (Theorem \ref{t:Yn}), we conclude that
$$
\frac{2}{Z_r} \sum_{\mathscr{L}_s \le y\le \mathscr{L}_r} \ee^{-V(y)} 
\to
0,
\qquad\hbox{\rm in $\P^*$-probability.}
$$

\noindent This yields (\ref{E(T(far))}).

We now turn to the {\bf proof of (\ref{E(T(close,bad))})}. Recall that $E_\omega^{(r)}[L_{T_\varnothing^+} (y)] = \frac{\pi_r(y)}{\pi_r(\varnothing)} = \ee^{U(\varnothing) -U(y)}= \ee^{U(\varnothing)}[\ee^{-V(y)} + \sum_{z\in \T: \, {\buildrel \leftarrow \over z} =y} \ee^{-V(z)}]$ for $y<\mathscr{L}_s$. By definition of $T_\varnothing^{\,\mathrm{(b)}}$ in (\ref{T(close,bad)}),
$$
E_\omega^{(r)}(T_\varnothing^{\,\mathrm{(b)}})
=
\ee^{U(\varnothing)} \sum_{y < \mathscr{L}_s} 
\Big( \ee^{-V(y)} + \sum_{z\in \T: \, {\buildrel \leftarrow \over z} =y} \ee^{-V(z)} \Big) \, 
{\bf 1}_{\{ y \; \mathrm{bad} \} } \, ,
$$

\noindent where, by ``$y$ bad", we mean $\min_{u\in \, [\! [ \varnothing, \, y]\! ]} \omega(u, \, {\buildrel \leftarrow \over u}) < (\log r)^{-6/\delta_1}$. So
$$
E_\omega^{(r)}(T_\varnothing^{\,\mathrm{(b)}})
\le
2\ee^{U(\varnothing)}\sum_{y \le \mathscr{L}_s} \ee^{-V(y)} {\bf 1}_{\{ y \; \mathrm{bad} \} } 
\le
2\ee^{U(\varnothing)}\sum_{y \le \mathscr{L}_r} \ee^{-V(y)} {\bf 1}_{\{ y \; \mathrm{bad} \} } \, .
$$

\noindent Let $B>0$ be a constant. Then
$$
E_\omega^{(r)}(T_\varnothing^{\,\mathrm{(b)}})
\le
2\ee^{U(\varnothing)} (\Sigma_{(\ref{Sigma1})} + \Sigma_{(\ref{Sigma2})}) \, ,
$$

\noindent where\footnote{For notational convenience, we treat $B(\log r)^2$ as an integer.}
\begin{eqnarray}
    \Sigma_{(\ref{Sigma1})}   
 &:=&\sum_{y \in \T: \; |y| \le B(\log r)^2} 
    \ee^{-V(y)} \, 
    {\bf 1}_{\{ y \; \mathrm{bad} \} }  \, ,
    \label{Sigma1}
    \\
    \Sigma_{(\ref{Sigma2})}   
 &:=&\sum_{y \le \mathscr{L}_r, \; |y| > B(\log r)^2} 
    \ee^{-V(y)} \, .
    \label{Sigma2}
\end{eqnarray}

\noindent In view of (\ref{Zn-proof6}), we have, for any $\varepsilon>0$,
\begin{equation}
    \lim_{B\to \infty} \, \limsup_{r\to \infty} \, 
    \P^* \Big\{ \Sigma_{(\ref{Sigma2})} \ge \varepsilon \log r\Big\}
    =
    0 \, .
    \label{Sigma2<}
\end{equation}

We now bound $\Sigma_{(\ref{Sigma1})}$. When $y$ is bad, $\max_{u\in \, [\! [ \varnothing, \, y]\! ]} \frac{1}{\omega(u, \, {\buildrel \leftarrow \over u})} > (\log r)^{6/\delta_1}$, which means $\max_{u\in \, [\! [ \varnothing, \, y]\! ]} \Lambda (u) > (\log r)^{6/\delta_1}-1$, with $\Lambda(u) := \sum_{z\in \T: \, {\buildrel \leftarrow \over z} =u} \ee^{-[V(z)-V(u)]}$ as in (\ref{Lambda}). Accordingly, writing $a(r) := (\log r)^{6/\delta_1}-1$ for brevity,
\begin{eqnarray*}
    \Sigma_{(\ref{Sigma1})}
 &\le&\sum_{y\in \T: \, |y| \le B(\log r)^2} 
    \ee^{-V(y)} \, 
    {\bf 1}_{\{ \max_{u\in \, [\! [ \varnothing, \, y]\! ]} \Lambda(u)> a(r) \} }
    \\
 &\le& 1+
    \sum_{k=1}^{B(\log r)^2} \sum_{y\in \T: \, |y| =k} 
    \ee^{-V(y)} \, 
    \Big( {\bf 1}_{\{ \max_{u\in \, [\! [ \varnothing, \, y[\![} \Lambda(u)> a(r) \} }
    + 
    {\bf 1}_{\{ \Lambda(y)> a(r) \} } \Big) \, ,
\end{eqnarray*}

\noindent the first term ``1" on the right-hand side resulting from $y= \varnothing$. We take expectation with respect to $\P$ on both sides. By Lemma \ref{l:many-to-one-appli1} and in its notation,
$$
\E(\Sigma_{(\ref{Sigma1})})
\le
1
+
\sum_{k=1}^{B(\log r)^2} \E\Big[ {\bf 1}_{\{ \max_{1\le i\le k} \widetilde{\Lambda}_{i-1} > a(r) \} } 
+ 
\P\Big( \sum_{x: \, |x|=1} \ee^{-V(x)} > a(r) \Big) \Big] \, .
$$

\noindent Since $\widetilde{\Lambda}_{i-1}$ (for $i\ge 1$) is distributed as $\widetilde{\Lambda}_0$, we have, for all $b>0$ and all $i\ge 1$, $\P(\widetilde{\Lambda}_{i-1} > b) \le c_8 \, b^{-\delta_1}$, where $\delta_1>0$ is the constant in (\ref{integrabilite_sous_Q}), and $c_8 := \E [ ( \widetilde{\Lambda}_0)^{\delta_1} ] = \E [ (\sum_{x: \, |x|=1} \ee^{-V(x)})^{1+\delta_1} ] $ which is finite according to (\ref{integrabilite_sous_Q}). On the other hand, (\ref{integrabilite_sous_Q}) also yields $\P(\sum_{x: \, |x|=1} \ee^{-V(x)} > b) \le c_8 \, b^{-(1+\delta_1)}$ (for $b>0$). Hence
$$
\E(\Sigma_{(\ref{Sigma1})})
\le
1
+
B(\log r)^2 \Big[ B(\log r)^2 \, c_8 \, a(r)^{-\delta_1} + c_8 \, a(r)^{-(1+\delta_1)} \Big] \, .
$$

\noindent This yields $\frac{\Sigma_{(\ref{Sigma1})}}{\log r} \to 0$ in $L^1(\P)$ and equivalently, in $L^1(\P^*)$, and a fortiori in $\P^*$-probability. Together with (\ref{Sigma2<}), and since $E_\omega^{(r)}(T_\varnothing^{\,\mathrm{(b)}}) \le 2 \ee^{U(\varnothing)} ( \Sigma_{(\ref{Sigma1})} + \Sigma_{(\ref{Sigma2})})$, we obtain $\frac{E_\omega^{(r)}(T_\varnothing^{\,\mathrm{(b)}})}{\log r} \to 0$ in $\P^*$-probability. Recalling that $E_\omega^{(r)} (T_\varnothing^+) = \frac{1}{\pi_r(\varnothing)}$ and that $(\log r)\pi_r(\varnothing)$ converges in $\P^*$-probability to a positive limit (Theorem \ref{t:Yn}), we deduce that
$$
\frac{E_\omega^{(r)} (T_\varnothing^{\,\mathrm{(b)}})}{E_\omega^{(r)} (T_\varnothing^+)}
\to 0,
\qquad \hbox{in $\P^*$-probability},
$$    

\noindent which is the desired conclusion in (\ref{E(T(close,bad))}). Lemma \ref{l:E(T)} is proved.\hfill$\Box$

\bigskip

\noindent {\it Proof of Lemma \ref{l:V(T)}.} Recall that $T_\varnothing^{\,\mathrm{(a)}} := \sum_{y < \mathscr{L}_s} L_{T_\varnothing^+} (y) \, {\bf 1}_{\{ y \; \mathrm{good}\} }$, where
$$
\{ y \; \mathrm{good}\}
:=
\Big\{ \min_{u\in \, [\! [ \varnothing, \, y]\! ]} \omega(u, \, {\buildrel \leftarrow \over u}) \ge (\log r)^{-6/\delta_1} \Big\} \, ,
$$

\noindent with $\delta_1>0$ denoting the constant in (\ref{integrabilite_sous_Q}). For any $y< \mathscr{L}_s$, we have $y<\mathscr{L}_r$, so $T_\varnothing^{\,\mathrm{(a)}}$ has the same distribution under $P_\omega^{(r)}$ and under $P_\omega$. In particular,
$$
\mathrm{Var}_\omega^{(r)} (T_\varnothing^{\,\mathrm{(a)}})
=
\mathrm{Var}_\omega (T_\varnothing^{\,\mathrm{(a)}})
=
\mathrm{Var}_\omega \Big( \sum_{y < \mathscr{L}_s} L_{T_\varnothing^+} (y) \, {\bf 1}_{\{ y \; \mathrm{good}\} } \Big) \, .
$$

\noindent Let $\overline{L}_k(x) := \sum_{i=1}^k {\bf 1}_{\{ X_{i-1} = {\buildrel \leftarrow \over x} , \, X_i =x\} } $ be edge local time as in (\ref{edge_local_time}). Then
$$
\sum_{y < \mathscr{L}_s} L_{T_\varnothing^+} (y) \, {\bf 1}_{\{ y \; \mathrm{good}\} }
=
\sum_{y < \mathscr{L}_s} \overline{L}_{T_\varnothing^+} (y) \, {\bf 1}_{\{ y \; \mathrm{good}\} }
+
\sum_{y \le \mathscr{L}_s, \, y\not= \varnothing} \overline{L}_{T_\varnothing^+} (y) \, {\bf 1}_{\{ {\buildrel \leftarrow \over y} \; \mathrm{good}\} } \, .
$$

\noindent By the elementary inequality $\mathrm{Var}(\xi_1 + \xi_2) \le 2[ \mathrm{Var}(\xi_1) + \mathrm{Var}(\xi_2)]$ (for random variables $\xi_1$ and $\xi_2$ having finite second moments), this leads to:
$$
\mathrm{Var}_\omega^{(r)} (T_\varnothing^{\,\mathrm{(a)}})
\le
2\, \mathrm{Var}_\omega \Big( \sum_{y < \mathscr{L}_s} \overline{L}_{T_\varnothing^+} (y) \, {\bf 1}_{\{ y \; \mathrm{good}\} } \Big)
+
2\, \mathrm{Var}_\omega \Big( \sum_{y \le \mathscr{L}_s, \, y\not= \varnothing} \overline{L}_{T_\varnothing^+} (y) \, {\bf 1}_{\{ {\buildrel \leftarrow \over y} \; \mathrm{good}\} } \Big) \, .
$$

\noindent We write
\begin{eqnarray*}
    \mathrm{Var}_\omega \Big( \sum_{y \le \mathscr{L}_s, \, y\not= \varnothing} \overline{L}_{T_\varnothing^+} (y) \, {\bf 1}_{\{ {\buildrel \leftarrow \over y} \; \mathrm{good}\} } \Big) 
 &=& \sum_{y \le \mathscr{L}_s, \, y\not= \varnothing} 
    \mathrm{Var}_\omega [\overline{L}_{T_\varnothing^+} (y)]
    \, {\bf 1}_{\{ {\buildrel \leftarrow \over y} \; \mathrm{good}\} } 
    +
    \\
 && \hskip-50pt +
    \sum_{y \not= z \le \mathscr{L}_s, \, y, \, z\not= \varnothing} 
    \mathrm{Cov}_\omega 
    [\overline{L}_{T_\varnothing^+} (y), \, \overline{L}_{T_\varnothing^+} (z) ] 
    \, {\bf 1}_{\{ {\buildrel \leftarrow \over y} \; \mathrm{good}\} }
    \, {\bf 1}_{\{ {\buildrel \leftarrow \over z} \; \mathrm{good}\} } \, , 
\end{eqnarray*}

\noindent and we have a similar expression for $\mathrm{Var}_\omega ( \sum_{y < \mathscr{L}_s} \overline{L}_{T_\varnothing^+} (y) \, {\bf 1}_{\{ y \; \mathrm{good}\} } )$.
Lemma \ref{l:V(T)} will be a straightforward consequence of the following inequalities: for some constants $c_9>0$ and $c_{10}>0$, and all $r\ge 2$,
\begin{eqnarray}
    \E \Big( \sum_{y \le \mathscr{L}_s} 
    E_\omega [\overline{L}_{T_\varnothing^+} (y)^2] \, {\bf 1}_{\{ {\buildrel \leftarrow \over y} \; \mathrm{good}\} } \Big)
 &\le& c_9\, s \, (\log r)^{\frac{6}{\delta_1} + 2} \, ,
    \label{second-moment:1st-technical-inequality}
    \\
    \E \Big( \sum_{y \not= z \le \mathscr{L}_s} 
   \big( \mathrm{Cov}_\omega 
    [\overline{L}_{T_\varnothing^+} (y), \, \overline{L}_{T_\varnothing^+} (z) ] \big)^+
    \, {\bf 1}_{\{ {\buildrel \leftarrow \over y} \; \mathrm{good}\} } \Big)
 &\le& c_{10} \, s\, (\log r)^{\frac{18}{\delta_1}+6} \, ,
    \label{second-moment:2nd-technical-inequality}
\end{eqnarray}

\noindent where $\delta_1>0$ is in (\ref{integrabilite_sous_Q}), and $(\mathrm{Cov}_\omega [\cdots ])^+$ denotes the positive part of $\mathrm{Cov}_\omega [\cdots]$.

So it remains to check inequalities (\ref{second-moment:1st-technical-inequality}) and (\ref{second-moment:2nd-technical-inequality}). We start with the {\bf proof of (\ref{second-moment:1st-technical-inequality})}. Recall from Lemma \ref{l:Cov} that
$$
E_\omega [\overline{L}_{T_\varnothing^+} (y)^2]
=
\omega(\varnothing, \, {\buildrel \leftarrow \over \varnothing})\, \ee^{-V(y)} \Big( 2 \sum_{z\in \, ]\!] \varnothing, \, y]\!]} \ee^{V(z)-V(y)} -1\Big)
\le
2\, \ee^{-V(y)} \sum_{z\in \, ]\!] \varnothing, \, y]\!]} \ee^{V(z)-V(y)} \, .
$$

\noindent For the sum on the right-hand side, we write (${\buildrel \Leftarrow \over y}$ denoting as before the parent of ${\buildrel \leftarrow \over y}$)
\begin{eqnarray*}
    \sum_{z\in \, ]\!] \varnothing, \, y]\!]} \ee^{V(z)-V(y)}
 &=&\ee^{-[V(y)-V({\buildrel \leftarrow \over y})]} \sum_{z\in \, ]\!] \varnothing, \, {\buildrel \leftarrow \over y}]\!]} \ee^{V(z)-V({\buildrel \leftarrow \over y})}
    +
    1
    \\
 &=&\frac{\omega({\buildrel \leftarrow \over y}, \, y)}{\omega({\buildrel \leftarrow \over y}, \, {\buildrel \Leftarrow \over y})} \sum_{z\in \, ]\!] \varnothing, \, {\buildrel \leftarrow \over y}]\!]} \ee^{V(z)-V({\buildrel \leftarrow \over y})}
    +
    1
    \\
 &\le&\frac{1}{\omega({\buildrel \leftarrow \over y}, \, {\buildrel \Leftarrow \over y})} \sum_{z\in \, ]\!] \varnothing, \, {\buildrel \leftarrow \over y}]\!]} \ee^{V(z)-V({\buildrel \leftarrow \over y})}
    +
    1\, .
\end{eqnarray*}

\noindent If $y \le \mathscr{L}_s$, then by definition, ${\buildrel \leftarrow \over y} < \mathscr{L}_s$, so that $\sum_{z\in \, ]\!] \varnothing, \, {\buildrel \leftarrow \over y}]\!]} \ee^{V(z)-V({\buildrel \leftarrow \over y})} \le s$. On the other hand, if ${\buildrel \leftarrow \over y}$ is good, then by definition, $\frac{1}{\omega({\buildrel \leftarrow \over y}, \, {\buildrel \Leftarrow \over y})} \le (\log r)^{6/\delta_1}$. Consequently, 
\begin{equation}
    E_\omega [\overline{L}_{T_\varnothing^+} (y)^2] \, {\bf 1}_{\{ {\buildrel \leftarrow \over y} \; \mathrm{good}\} } \, {\bf 1}_{\{ y \le \mathscr{L}_s \} } 
    \le 
    2\, [s \, (\log r)^{\frac{6}{\delta_1}} +1] \, \ee^{-V(y)} \, .
    \label{second-moment:3rd-technical-inequality}
\end{equation}

\noindent Since $\E(\, \sum_{y \le \mathscr{L}_s} \ee^{-V(y)}) \le c_2 \, (\log s)^2$ (see (\ref{E(Yn)<})), this yields (\ref{second-moment:1st-technical-inequality}).

We now turn to the {\bf proof of (\ref{second-moment:2nd-technical-inequality})}. Consider a pair $y \not= z \le \mathscr{L}_s$. By Lemma \ref{l:Cov},
$$
    \mathrm{Cov}_\omega 
    [\overline{L}_{T_\varnothing^+} (y), \, \overline{L}_{T_\varnothing^+} (z) ]
    \le 
    2\, \ee^{-[V(y)-V(y\wedge z)]-[V(z)-V(y\wedge z)]}\,
    E_\omega [\overline{L}_{T_\varnothing^+}(y\wedge z)^2] \, .
$$

\noindent Hence, writing $\mathrm{LHS}_{(\ref{second-moment:2nd-technical-inequality})} := \sum_{y \not= z \le \mathscr{L}_s} (\mathrm{Cov}_\omega [\overline{L}_{T_\varnothing^+} (y), \, \overline{L}_{T_\varnothing^+} (z) ])^+ \, {\bf 1}_{\{ {\buildrel \leftarrow \over y} \; \mathrm{good}\} }$, we have
$$
\mathrm{LHS}_{(\ref{second-moment:2nd-technical-inequality})}
\le
2\sum_{u < \mathscr{L}_s}
\sum_{y\not= z \le \mathscr{L}_s: \, y\wedge z =u} 
\ee^{-[V(y)-V(u)]-[V(z)-V(u)]}\, E_\omega [\overline{L}_{T_\varnothing^+}(u)^2] 
\, {\bf 1}_{\{ u \; \mathrm{good}\} }\, .
$$

\noindent Observe that ${\bf 1}_{\{ u \; \mathrm{good}\} } \le {\bf 1}_{\{ {\buildrel \leftarrow \over u} \; \mathrm{good}\} }$, so by (\ref{second-moment:3rd-technical-inequality}), $E_\omega [\overline{L}_{T_\varnothing^+} (u)^2] \, {\bf 1}_{\{ u \; \mathrm{good}\} } \le 2[s \, (\log r)^{\frac{6}{\delta_1}} +1] \, \ee^{-V(u)}$. Consequently,
\begin{eqnarray}
    \mathrm{LHS}_{(\ref{second-moment:2nd-technical-inequality})}
 &\le& 4[s \, (\log r)^{\frac{6}{\delta_1}} +1]
    \sum_{u < \mathscr{L}_s}
    \ee^{-V(u)} \, {\bf 1}_{\{ u \; \mathrm{good}\} } \times
    \nonumber
    \\
 &&\qquad \times
    \sum_{y\not= z\le\mathscr{L}_s: \, y\wedge z =u}
    \ee^{-[V(y)-V(u)]-[V(z)-V(u)]}\, .
    \label{LHS<}
\end{eqnarray}

\noindent Let us consider the double sum $\sum_{y\not= z\le\mathscr{L}_s: \, y\wedge z =u} \ee^{-[V(y)-V(u)]-[V(z)-V(u)]}$ on the right-hand side. Write $k=k(u) := |u|$ for brevity. Then
\begin{eqnarray*}
 &&\sum_{y\not= z\le\mathscr{L}_s: \, y\wedge z =u}
    \ee^{-[V(y)-V(u)]-[V(z)-V(u)]}
    \\
 &=&\sum_{a\not= b, \, {\buildrel \leftarrow \over a} = u = {\buildrel \leftarrow \over b}}
    \ee^{-[V(a)-V(u)]-[V(b)-V(u)]}
    \sum_{y, \, z\le \mathscr{L}_s: \, y_{k+1}=a, \, z_{k+1}=b}
    \ee^{-[V(y)-V(a)]-[V(z)-V(b)]} \, .
\end{eqnarray*}

\noindent Observe that if $y \le \mathscr{L}_s$ is such that $|y| \ge k+1$ and $y_{k+1}=a$, then by definition of $\mathscr{L}_s$ in (\ref{gamma}), $\sum_{w\in \, ]\! ] a, \, v]\! ]} \ee^{V(w)- V(v)} \le s$, $\forall v\in \, ]\!] a, \, y[\![\,$; so writing $y=a\widetilde{y}$ (the concatenation of $a$ and $\widetilde{y}$), then $\widetilde{y}$ as a vertex of the subtree rooted at $a$ satisfies $\widetilde{y} \le\mathscr{L}_s(a)$, where $\mathscr{L}_s(a)$ is defined as $\mathscr{L}_s$, but associated with the subtree rooted at vertex $a$. Accordingly, with $\mathscr{F}_{k+1}$ denoting the $\sigma$-field generated by $(V(x), \; |x| \le k+1)$, we have that on the set $\{ | u| =k\}$, 
\begin{eqnarray*}
 &&\E\Big( \, \sum_{y\not= z\le \mathscr{L}_s: \, y\wedge z =u}
    \ee^{-[V(y)-V(u)]-[V(z)-V(u)]} \, \Big| \, \mathscr{F}_{k+1} \Big)
    \\
 &\le& \sum_{a\not= b, \, {\buildrel \leftarrow \over a} = u= {\buildrel \leftarrow \over b}}
    \ee^{-[V(a)-V(u)]-[V(b)-V(u)]} \,
    [ \E (Y_s)]^2 \, .
\end{eqnarray*}

\noindent If $u$ is good, then by definition, $\omega(u, \, {\buildrel \leftarrow \over u}) \ge (\log r)^{-\frac{6}{\delta_1}}$; in particular,
$$
\sum_{a\in \T: \, {\buildrel \leftarrow \over a} = u} \ee^{-[V(a)-V(u)]}
=
\frac{1}{\omega(u, \, {\buildrel \leftarrow \over u})} -1
\le
\frac{1}{\omega(u, \, {\buildrel \leftarrow \over u})}
\le
(\log r)^{\frac{6}{\delta_1}} \, ,
$$

\noindent which implies that $\sum_{a\not= b, \, {\buildrel \leftarrow \over a} = u= {\buildrel \leftarrow \over b}} \ee^{-[V(a)-V(u)]-[V(b)-V(u)]} \le (\log r)^{12/\delta_1}$. Hence on the set $\{ | u| =k\}$,  
$$
\E\Big( \, \sum_{y\not= z\le \mathscr{L}_s: \, y\wedge z =u}
\ee^{-[V(y)-V(u)]-[V(z)-V(u)]} \, \Big| \, \mathscr{F}_{k+1} \Big)
\, {\bf 1}_{\{ u \; \mathrm{good}\} }
\le
(\log r)^{\frac{12}{\delta_1}}\, [ \E (Y_s)]^2 \, .
$$

\noindent Going back to (\ref{LHS<}), this yields
$$
\E(\mathrm{LHS}_{(\ref{second-moment:2nd-technical-inequality})})
\le
\E \Big( 4 \, [s \, (\log r)^{\frac{6}{\delta_1}} +1] (\log r)^{\frac{12}{\delta_1}} \, [ \E (Y_s)]^2\, \sum_{u < \mathscr{L}_s} \ee^{-V(u)}\Big) \, .
$$

\noindent Since $\sum_{u < \mathscr{L}_s} \ee^{-V(u)} \le \sum_{u \le \mathscr{L}_s} \ee^{-V(u)} = Y_s$, we obtain: $\E(\mathrm{LHS}_{(\ref{second-moment:2nd-technical-inequality})}) \le 4 \, [s \, (\log r)^{\frac{6}{\delta_1}} +1] (\log r)^{\frac{12}{\delta_1}} \, [ \E (Y_s)]^3$. In view of (\ref{E(Yn)<}), this yields (\ref{second-moment:2nd-technical-inequality}), and completes the proof of Lemma \ref{l:V(T)}.\hfill$\Box$

%
%
%

\section{Biased walks: proof of Theorem \ref{t:ligne-arret}} 
 \label{s:proof-t:ligne-arret} 

$\phantom{aob}$Recall from (\ref{pf:l:main}) that
$$
    P_\omega \Big( \bigcup_{i=1}^n\{ X_i \in \mathscr{L}_r\} \Big)
    \le
    \frac{E_\omega (L_n(\varnothing)+1)}{r} \sum_{x\in \mathscr{L}_r}\ee^{-V(x)} .
$$

\noindent By Lemma \ref{l:sum(xinLn)}, $(\log r) \sum_{x\in \mathscr{L}_r}\ee^{-V(x)}$ is tight under $\P^*$. Theorem \ref{t:ligne-arret} follows from \eqref{localtimemoment} of Proposition \ref{p:x=0}.\hfill$\Box$ 


\section{Biased walks: proof of Proposition \ref{p:local}} 
 \label{s:localproba}
 
$\phantom{aob}$We begin with a general fact for reversible Markov chains. The fact is well known. For a simple proof for finite chains, see Saloff-Coste (\cite{Saloff97}, Lemma 1.3.3~(1), page 323), applied to $P^2$.

\medskip

\begin{fact}
 \label{monotone}
 
 Let $P$ be the transition probability of a reversible Markov chain taking values in a countable space $E$. Then for any $x \in E$, the sequence $k \to P^{2 k}(x, \, x)$ is non-increasing. 

\end{fact}

\medskip






We prepare for the proof of Proposition \ref{p:local}. Let $P_\omega^{(r)}$ be, as before, the quenched probability with a reflecting barrier at $\mathscr{L}_r$, and $E_\omega^{(r)}$ the corresponding expectation.  

\medskip

\begin{lemma}
 \label{l:pn} 
 
 Let $\gamma \in \r$, and let $r=r(n) := \frac{n}{(\log n)^\gamma}$. Then 
 $$
 (\log n)^{2- \gamma} \,  \sup_{ B \in \sigma\{ X_1, \, \cdots \, , X_n\}} | P_\omega^{(r)} (B) - P_\omega(B) |
 $$
 is tight under $\P^*$.

\end{lemma}

\medskip
 
{\noindent\bf Proof.} For $B \in \sigma\{ X_1, \, \cdots \, , X_n\}$,  
$$ 
| P_\omega^{(r)} (B) - P_\omega(B) | 
\le
P_\omega \Big( \bigcup_{i=1}^n\{ X_i \in \mathscr{L}_r\} \Big) \, ,
$$ 

\noindent which is bounded by $\frac{E_\omega (L_n(\varnothing)+1)}{r} \sum_{x\in \mathscr{L}_r}\ee^{-V(x)}$ (see (\ref{pf:l:main})). We conclude by means of Lemma \ref{l:sum(xinLn)} and \eqref{localtimemoment} of Proposition \ref{p:x=0}.\hfill$\Box$

\bigskip

We are now ready to prove Proposition \ref{p:local}.

\bigskip

\noindent {\it Proof of Proposition \ref{p:local}.} We choose $r:=n$ so that we are entitled to apply Lemma \ref{l:pn} (with $\gamma =0$). We claim that for any $a_n\to \infty$ satisfying $\lim_{n \to \infty} \frac{\log a_n}{\log n}=0$,
\begin{equation} 
    \max_{k\; \mathrm{even}: \, \frac{n}{a_n} \le k \le n}   
    \Big| ( \log n)  P_\omega^{(r)} (X_k= \varnothing) - \frac{\sigma^2}{2 D_\infty} \ee^{- U(\varnothing)} \Big| 
    \to 
    0,  \qquad\hbox{in $\P^*$-probability} .
    \label{localn1}
\end{equation}

\noindent By Lemma \ref{l:pn}, \eqref{localn1} will imply Proposition \ref{p:local}.

Let $m:=m_n$ be the smallest even number such that $m \ge \frac{n}{a_n}$. Clearly $\frac{\log m}{\log n} \to 1$.  Using the trivial upper bound $\frac{L_m(\varnothing)}{m} \le 1$, we deduce from part \eqref{localtimemoment} of Proposition \ref{p:x=0} and Lemma \ref{l:pn} that for $n \to\infty$, 
\begin{equation}
    \label{localtimemoment2}  
    E_\omega^{(r)} \Big( \frac{L_m(\varnothing)}{\frac{m}{\log m}} \Big) 
    \to  
    \frac{\sigma^2}{4 D_\infty} \, \ee^{-U(\varnothing)},  \qquad\hbox{in $\P^*$-probability} \, ,
\end{equation}

By Fact \ref{monotone}, $i \mapsto P_\omega^{(r)}(X_{2i} = \varnothing)$ is non-increasing, so $E_\omega^{(r)} ( \frac{L_m(\varnothing)}{\frac{m}{\log m}} ) = \frac{\log m}{m} \sum_{i=1}^m P_\omega^{(r)} (X_i= \varnothing)\ge \frac12 (\log m) P_\omega^{(r)} ( X_m= \varnothing )$, the factor $\frac12$ coming from the fact we sum over even numbers $i\in [1, \, m]$. Combined with \eqref{localtimemoment2}, we see that
\begin{equation}
    \max_{k \; \mathrm{even}: \; \frac{n}{a_n} \le k \le n} 
    ( \log n)  P_\omega^{(r)} ( X_k= \varnothing )   
    =
    (\log n)  P_\omega^{(r)} ( X_m= \varnothing )   
    \le  
    \frac{\sigma^2 +o_{\P^*}(1)}{2 D_\infty} \ee^{- U(\varnothing)}, 
    \label{localupper}
\end{equation}
 
\noindent where $o_{\P^*}(1)$ denotes a quantity which goes to $0$ in $\P^*$-probability as $n \to \infty$.
 
To obtain the lower bound for $P_\omega^{(r)} ( X_k= \varnothing )$,  we consider the Markov chain $(X_{2i}, \, i\ge 0)$ under $P_\omega^{(r)}$, starting from $X_0:=\varnothing$. This chain takes values in $E_r:= \{ x\in \T:  x \le \mathscr{L}_r, \, |x| \;\mathrm{even}\}$, with $\pi_r(E_r)=\frac12$ due  to periodicity. In other words, $2 \pi_r(x)$ for $x\in E_r$, is the invariant probability measure of $(X_{2i}, \, i\ge 0 )$. By Fact \ref{monotone}, we see that for integer $i\ge 0$, $P_\omega^{(r)} (X_{2i}= \varnothing) \ge 2 \, \pi_r(\varnothing)$. In particular, for $k:=2\lfloor \frac{n}{2} \rfloor$, $P_\omega^{(r)} ( X_k= \varnothing ) \ge 2 \pi_r(\varnothing)= \frac{2}{Z_r} \ee^{-U(\varnothing)}$. As such,   \eqref{localn1} follows from Theorem \ref{t:Yn} and \eqref{localupper}.\hfill$\Box$
  
\medskip

\begin{remark}

By definition, $\frac{1}{\pi_r(\varnothing)} = Z_r \, \ee^{U(\varnothing)}$, which is $\frac{4+o_{\P^*}(1)}{\sigma^2} D_\infty \ee^{U(\varnothing)} \log n$ according to Theorem \ref{t:Yn}, where $o_{\P^*}(1)\to 0$ in $\P^*$-probability as $n \to \infty$. So \eqref{localn1} can also be stated as follows: For any $a_n\to \infty$ such that $\lim_{n \to \infty} \frac{\log a_n}{\log n}=0$, uniformly in even integers $k\in [\frac{n}{a_n}, \, n]$, 
\begin{equation}
    P_\omega^{(n)} (X_k= \varnothing) 
    =
    (2+o_{\P^*}(1))  \pi_n(\varnothing)\, . 
    \label{localn2}
\end{equation}
  
\noindent This will be useful in the proof of Theorem \ref{t:cv} in Section \ref{s:maintheorem}. 

\end{remark}

\section{Biased walks: proofs of Lemma \ref{l:pinpim}, Theorem \ref{t:cv} and Corollary \ref{c:cv}}
 \label{s:maintheorem}
 
{\noindent\it Proof of Lemma \ref{l:pinpim}.} For $0<u\le r$ and $x \in \T \cup\{{\buildrel \leftarrow \over \varnothing}\}$,  
\begin{eqnarray*}
    | \pi_r(x)- \pi_u(x) |
 &=& \Big| \frac{Z_r \pi_r(x) - Z_u \pi_u(x)}{Z_u} - Z_r \pi_r(x) \big( \frac{1}{Z_u} - \frac{1}{Z_r} \big) \Big| 
    \\
 &\le& \frac{Z_r \pi_r(x) - Z_u \pi_u(x)}{Z_u} +  Z_r \pi_r(x) \big( \frac{1}{Z_u} - \frac{1}{Z_r} \big),
\end{eqnarray*}

\noindent by using the facts that $Z_r \pi_r(x) \ge  Z_u \pi_u(x)$ and $Z_r \ge Z_u$. By taking the summation on $x$,   
$$ 
2\, \dtv (\pi_u, \, \pi_r) 
\le  
\frac{Z_r-Z_u}{Z_u} + Z_r \big( \frac{1}{Z_u} - \frac{1}{Z_r} \big) 
=
2 \frac{Z_r - Z_u}{Z_u} \, ,
$$

\noindent which is bounded by $2 \frac{Z_r - Z_{r/(\log r)^a}}{Z_{r/(\log r)^a}}$ if $u\in [\frac{r}{(\log r)^a}, \, r]$. By Theorem \ref{t:Yn}, $\frac{1}{\log r} \big( Z_r - Z_{r/(\log r)^a}\big)\to 0$ in ${\bf P^*}$-probability, from which Lemma \ref{l:pinpim} follows.\hfill$\Box$
 
\medskip 

\noindent {\it Proof of Theorem \ref{t:cv}.} 

{\bf (i) Case $\kappa=1$.} We prove the following  stronger statement: Fix $0< c<1$. As $n\to \infty$,  
\begin{equation}
    \max_{\lfloor cn\rfloor \le m \le n}  
    \sup_{A \subset \T \cup\{ {\buildrel \leftarrow \over \varnothing} \}}
    | P_\omega (X_m \in A) - \widetilde{\pi}_m (A) |   
    \to 
    0, \qquad  \hbox{in $\P^*$-probability} .
    \label{xnpi-new}
\end{equation} 
 
\noindent The fact that \eqref{xnpi-new} holds uniformly in $m$ will be useful in the proof for the case $\kappa\ge2$. 

 
In view of Lemma \ref{l:pn}, it suffices to prove \eqref{xnpi-new} for $P_\omega^{(r)}$ in lieu of $P_\omega$, with $r:= n$. 

Let $n$ be large and put $b_n:= \lfloor \frac{n}{(\log n)^2} \rfloor$. For $cn\le m \le n$ (we treat $cn$ as an integer) and $A \subset \T \cup\{{\buildrel \leftarrow \over \varnothing}\}$, we have
\begin{equation}
    0
    \le 
    P_\omega^{(r)} (X_m \in A) - \sum_{k=b_n}^m P_\omega^{(r)} (X_m \in A, \, \mathfrak{g}_m=k)
    \le 
    P_\omega^{(r)} \{ L_{cn}(\varnothing) < b_n \} \, ,
    \label{pf:t:main-1}
\end{equation}

\noindent  where $\mathfrak{g}_m:= \max\{i \le m: X_i= \varnothing\}$ is the last return time to $\varnothing$ before $m$ and the second inequality follows from the fact that $\{ \mathfrak{g}_m \le b_n\} \subset \{ L_m(\varnothing) \le b_n \} \subset \{ L_{c n} (\varnothing) \le b_n\}$.  

For any $\varepsilon>0$, we have $b_n < \varepsilon \frac{c n}{\log (c n)}$ for sufficiently large $n$; so
$$ 
P_\omega \{ L_{c n}(\varnothing) < b_n\}
\le  
P_\omega \Big( \, \Big| \frac{L_{cn}(\varnothing)}{\frac{c n}{\log (c n)}} - \frac{\sigma^2}{4 D_\infty} \, \ee^{-U(\varnothing)}\Big| > \varepsilon \Big) 
+
{\bf 1}_{\{ \frac{\sigma^2}{4 D_\infty} \, \ee^{-U(\varnothing)} \le 2\varepsilon\}}.
$$

\noindent Applying Proposition \ref{p:x=0}, and since $\varepsilon>0$ can be as small as possible, we see that $P_\omega \{ L_{c n}(\varnothing) < b_n\} \to 0$ in $\P^*$-probability. A fortiori, $P_\omega^{(r)} \{ L_{c n}(\varnothing) < b_n\} \to 0$ in $\P^*$-probability. Going back to \eqref{pf:t:main-1}, we obtain:
\begin{equation}
    P_\omega^{(r)} (X_m \in A)   
    -
    \sum_{k =b_n}^m P_\omega^{(r)} (X_m \in A, \, \mathfrak{g}_m=k)  
    \to 
    0, \qquad  \hbox{in $\P^*$-probability} \, ,
    \label{lnsmall}
\end{equation}
  
\noindent uniformly in $A\subset \T \cup\{ {\buildrel \leftarrow \over \varnothing} \}$ and in $m\in [cn, \, n]\cap \z$.
 
Let us deal with the sum on the left-hand side of \eqref{lnsmall}. By the Markov property at time $k$,
$$ 
P_\omega^{(r)} (X_m \in A, \, \mathfrak{g}_m=k)
=
P_\omega^{(r)} (X_k=\varnothing) \, P_\omega^{(r)} (X_{m-k} \in A, \, m-k < T^+_\varnothing) \, ,
$$ 

\noindent where $T^+_\varnothing$ denotes, as before, the first return time to $\varnothing$.     
  
By \eqref{localn2}, uniformly in even numbers $k\in [b_n, \, n]$, $P_\omega^{(r)} (X_k= \varnothing) =  (2  +o_{\P^*}(1)) \pi_r(\varnothing)$. It follows that uniformly in $A$ and in $m$, 
$$  
P_\omega^{(r)} (X_m \in A)
-
2 \pi_r(\varnothing) \sum_{k =b_n, \; k\; \mathrm{even}}^m P_\omega^{(r)} (X_{m-k} \in A , \, m-k < T^+_\varnothing)
\to 
0 \, ,
$$   
   
\noindent in $\P^*$-probability. If $m$ is even, so is $m-k$, then we can restrict $A$ to $A \cap \T^{(\mathrm{even})}$. A similar restriction holds if $m$ is odd. Define 
$$    
A_m
:=
\begin{cases}
 A\cap \T^{(\mathrm{even})} \, , & \hbox{\rm if $m$ is even} \, ,
 \\ 
 A\cap (\T^{(\mathrm{odd})} \cup\{{\buildrel \leftarrow \over \varnothing}\}) \, , &\hbox{\rm if $m$ is odd}\, .    \end{cases}
$$

\noindent We have (with $o_{\P^*}(1)$ denoting an expression tending to $0$ in $\P^*$-probability, uniformly in $A$ and in $m$) 
\begin{eqnarray}  
    P_\omega^{(r)} (X_m \in A)
 &=& 2\pi_r(\varnothing) \, \sum_{k =b_n }^m P_\omega^{(r)} (X_{m-k} \in A_m , \, m-k < T^+_\varnothing)
    +
    o_{\P^*}(1)
    \nonumber
    \\
 &=&2\pi_r(\varnothing) \, \sum_{i=0}^{m-b_n} P_\omega^{(r)} (X_i \in A_m , \, i < T^+_\varnothing)
    +
    o_{\P^*}(1) \, ,
    \label{eq1}
\end{eqnarray}   

\noindent which implies that  
\begin{eqnarray}
    P_\omega^{(r)} (X_m \in A)
 &\le& 2\pi_r(\varnothing) \,  E_\omega^{(r)} \Big[ \sum_{i=0}^{T_\varnothing -1} {\bf 1}_{\{ X_i \in A_m\} } \Big] 
    +
    o_{\P^*}(1) 
    \label{eq2}
    \\ 
 &=&2 \, \pi_r (A_m) + o_{\P^*}(1) \, , 
    \nonumber
\end{eqnarray} 
  
\noindent by the fact that $\pi_r(\varnothing) = \frac{1}{E_\omega^{(r)} (T_\varnothing^+)}$. By Lemma \ref{l:pinpim}, $\pi_r (A_m) = \pi_m (A_m) + o_{\P^*}(1)$. Since $2 \pi_m (A_m) =  \widetilde \pi_m (A)$, we obtain that $P_\omega^{(r)} (X_m \in A) \le \widetilde{\pi}_m (A ) + o_{\P^*}(1)$. 
  
To get \eqref{xnpi-new}, it remains to check that $P_\omega^{(r)} (X_m \in A) \ge \widetilde{\pi}_m (A ) + o_{\P^*}(1)$, which will be done if we are able to reverse the inequality in  \eqref{eq2}. By \eqref{eq1} and tightness of $(\log n) \pi_r(\varnothing)$, it suffices to prove that 
\begin{equation} 
    \frac{1}{\log n} \sum_{i = m - b_n+1}^\infty 
    P_\omega^{(r)} (i < T^+_\varnothing)
    \to 
    0, \qquad\hbox{in $\P^*$-probability} .
    \label{eq3}
\end{equation}  

Of course, $\sum_{i = m- b_n+1}^\infty P_\omega^{(r)}(i < T^+_\varnothing) = E_\omega^{(r)} [ (T^+_\varnothing - (m-b_n+1))^+ ]$. By (\ref{T=Ta+Tb+Tc}) and in its notation (with $s:= \frac{r}{(\log r)^\theta}$ and $\theta\ge 0$), $T_\varnothing^+ \le T_\varnothing^{\,\mathrm{(a)}} + T_\varnothing^{\,\mathrm{(b)}} + T_\varnothing^{\,\mathrm{(c)}} +1$; so $\sum_{i = m- b_n+1}^\infty P_\omega^{(r)}(i < T^+_\varnothing) \le E_\omega^{(r)} [ (T^{(a)}_\varnothing - (m-b_n))^+ ] + E_\omega^{(r)} [T^{(b)}_\varnothing] + E_\omega^{(r)} [T^{(c)}_\varnothing]$. Lemma \ref{l:E(T)} entails that $E_\omega^{(r)} [T^{(b)}_\varnothing] + E_\omega^{(r)} [T^{(c)}_\varnothing] = (\log n) \times o_{\P^*}(1)$. On the other hand,
$$
E_\omega^{(r)} [(T^{(a)}_\varnothing - (m-b_n))^+]
\le 
E_\omega^{(r)} [T^{(a)}_\varnothing \, {\bf 1}_{\{ T^{(a)}_\varnothing \ge m-b_n\} }]
\le 
\frac{1}{m-b_n}\, E_\omega^{(r)} [(T^{(a)}_\varnothing)^2] \, .
$$ 
 
\noindent By Lemma \ref{l:V(T)} and tightness of $\frac{1}{\log n} E_\omega^{(r)} [T^{(a)}_\varnothing]$, we take a large parameter $\theta$ such that $\frac{18}{\delta_1}+6-\theta < 1$ and arrive at \eqref{eq3}. This completes the proof of \eqref{xnpi-new}.   

\medskip
  
{\bf (ii) Case $\kappa\ge 2$.} We only check the case $\kappa=2$ because the general case can be proved exactly in the same way. Without loss of generality, we take $t_2=1$ and $t_1=s\in (0, \, 1)$. For brevity, we treat $sn$ as an integer. It suffices to prove that, for $n\to \infty$,  
\begin{equation}
    \sup_{ A_1, \; A_2 \subset\T \cup\{{\buildrel \leftarrow \over \varnothing}\}}     
    | P_\omega (X_{sn} \in A_1, \, X_n \in A_2) - \widetilde{\pi}_{s n} (A_1) \, \widetilde{\pi}_n(A_2) |
    \to 
    0, \qquad  \hbox{in $\P^*$-probability} .
    \label{xnpi2}
\end{equation} 

Fix $t \in (s, \, 1)$. Let $\mathfrak{d}_{sn}:= \min\{ i>sn: \, X_i= \varnothing\}$ and $B_n:= \{ \mathfrak{d}_{s n } \le t n \} = \{ L_{t n} (\varnothing) > L_{sn} (\varnothing)\}$. By \eqref{localtimeproba} of Proposition \ref{p:x=0},
\begin{equation}
    P_\omega (B_n^c)
    \to 
    0, \qquad  \hbox{in $\P^*$-probability} .
    \label{bnc}
\end{equation}

\noindent Hence $P_\omega (X_{sn} \in A_1, \, X_n \in A_2) = P_\omega (X_{sn} \in A_1, \, X_n \in A_2, \, B_n ) + o_{\P^*}(1)$, where $o_{\P^*}(1)$ denotes an expression converging to $0$ in $\P^*$-probability uniformly in $A_1$, $A_2 \subset\T \cup\{{\buildrel \leftarrow \over \varnothing}\}$. Applying the strong Markov property at $\mathfrak{d}_{sn}$, this gives
$$
P_\omega (X_{sn} \in A_1, \, X_n \in A_2)
=
\sum_{k=sn+1}^{tn} P_\omega (X_{sn} \in A_1, \, \mathfrak{d}_{s n} =k) \, P_\omega (X_{n- k} \in A_2) 
+
o_{\P^*}(1) \, ,
$$

\noindent which is $\sum_{k= sn+1}^{t n} P_\omega (X_{sn} \in A_1, \, \mathfrak{d}_{s n} =k) \, \widetilde{\pi}_{n-k} (A_2) + o_{\P^*}(1)$ by \eqref{xnpi-new}. 

For even numbers $k\in (s n, \, tn]$, $n$ and $n-k$ have the same parity and $\dtv (\widetilde{\pi}_{n-k}, \, \widetilde{\pi}_n) \le 2 \dtv (\pi_{n-k},\, \pi_n)$. So by Lemma \ref{l:pinpim}, $\dtv (\widetilde{\pi}_{n-k}, \, \widetilde{\pi}_n) \to 0$ in $\P^*$-probability, uniformly in even numbers $k\in [s n, \, tn]$. As such,
\begin{eqnarray*} 
    P_\omega (X_{sn} \in A_1, \, X_n \in A_2)
 &=& \sum_{k=sn+1}^{tn} P_\omega (X_{sn} \in A_1, \, \mathfrak{d}_{sn} =k) \,  \widetilde{\pi}_{n} (A_2) 
    + 
    o_{\P^*}(1)  
    \\
 &=& P_\omega (X_{sn} \in A_1) \, \widetilde{\pi}_{n} (A_2) 
    + 
    o_{\P^*}(1) \, ,
\end{eqnarray*}

\noindent by means of \eqref{bnc}.  Applying the  already proved case $\kappa=1$ of Theorem \ref{t:cv} to $sn$, we get that $P_\omega (X_{sn} \in A_1) = \widetilde{\pi}_{s n} (A_1) + o_{\P^*}(1)$, which yields \eqref{xnpi2} and completes the proof of Theorem \ref{t:cv}.\hfill$\Box$
 
\bigskip

\noindent {\it Proof of Corollary \ref{c:cv}.} We only prove the case $\kappa=1$ and $t_\kappa=1$. The general case can be handled exactly in the same way. 

By Lemma \ref{l:pinpim}, $\dtv (\pi_n(\cdot), \, \pi_{n-1}(\cdot)) \to 0$ in $\P^*$-probability, from which follows that
\begin{equation} 
    \dtv \Big( \frac12 ( \widetilde{\pi}_n (\cdot) + \widetilde{\pi}_{n-1}(\cdot)), \, \pi_n(\cdot) \Big) 
    \to
    0,
    \qquad \hbox{in $\P^*$-probability}\, .
    \label{last5}
\end{equation}

Applying Theorem \ref{t:cv} (case $\kappa=1$) to $n$ and $n-1$, we get from \eqref{last5} that  
\begin{equation}
    \sup_{A \subset \T \cup\{{\buildrel \leftarrow \over \varnothing}\} }
    | P_\omega (X_n \in A) + P_\omega (X_{n-1} \in A) -  2 \pi_n (A) |
    \to 
    0, \qquad  \hbox{in $\P^*$-probability}\, .
    \label{xnpi-ter}
\end{equation} 
 
\noindent Let $B>b>0$ be constants and let $n$ be large. We treat $b(\log n)^2$ and $B (\log n)^2$ as integers for brevity. By \eqref{xnpi-ter} (with $o_{\P^*}(1)$ denoting an expression converging to $0$ in $\P^*$-probability, uniformly in $B>b>0$)
\begin{eqnarray}
    I_{\eqref{pnbB}}
 &:=& P_\omega  \Big( b\le \frac{|X_n|}{(\log n)^2} \le B \Big)
    +
    P_\omega  \Big( b\le \frac{|X_{n-1}|}{(\log n)^2} \le B \Big)
    \nonumber
    \\
 &=& 2 \sum_{b (\log n)^2 \le |x| \le B (\log n)^2} \pi_n(x) 
    + o_{\P^*}(1)\, .
    \label{pnbB}
\end{eqnarray}
 
\noindent By definition of $\pi_n$ in \eqref{pi_n} (and the fact that $Z_n = 2Y_n$ as in Lemma \ref{l:Zn=2Yn}),
\begin{eqnarray*}
    2 \sum_{b(\log n)^2 \le |x| \le B(\log n)^2} \pi_n(x)     
 &=& \frac{1}{Y_n} \sum_{k=b(\log n)^2}^{B(\log n)^2} \,
    \sum_{|x|=k}
    \Big( \, {\bf 1}_{\{ x< \mathscr{L}_n\} } \, \ee^{-U(x)}
    +
    {\bf 1}_{\{ x\in \mathscr{L}_n\} } \, \ee^{-V(x)} \, \Big)
    \\
 &=& \frac{2}{Y_n}\Big(  
    \sum_{k=b(\log n)^2}^{B(\log n)^2} \, 
    \sum_{|x|=k} 
    {\bf 1}_{\{ x< \mathscr{L}_n\} } \, \ee^{-V(x)} 
    +
    \Delta_n \Big), 
\end{eqnarray*}

\noindent where $|\Delta_n| \le \sum_{|x|=b(\log n)^2} \ee^{-V (x)} + \sum_{|x|=B(\log n)^2+1} \ee^{-V (x)} + \sum_{x\in \mathscr{L}_n} \ee^{-V(x)} = W_{b(\log n)^2} + W_{B(\log n)^2+1} + \sum_{x\in \mathscr{L}_n} \ee^{-V(x)}$, where $(W_i)$ is the additive martingale in \eqref{Wn}. Since $W_i\to 0$ (for $i\to \infty$) $\P^*$-a.s.\ (see \eqref{W->0}), and $\sum_{x\in \mathscr{L}_n} \ee^{-V(x)} \to 0$ in $\P^*$-probability (Lemma \ref{l:sum(xinLn)}), we have $\Delta_n \to 0$ in $\P^*$-probability. On the other hand, $\frac{Y_n}{\log n} \to \frac{2}{\sigma^2} D_\infty$ in $\P^*$-probability (Theorem \ref{t:Yn}). Consequently,
$$
I_{\eqref{pnbB}}
=
\frac{\sigma^2}{D_\infty\, \log n} \sum_{k=b(\log n)^2}^{B(\log n)^2} \, \sum_{|x|=k} {\bf 1}_{\{ x< \mathscr{L}_n\} } \, \ee^{-V(x)} 
+
o_{\P^*}(1) \, .
$$

\noindent By an obvious analogue of \eqref{pf-Yn-ub} and \eqref{pf-Yn-low}, 
$$
\sum_{k=b(\log n)^2}^{B(\log n)^2} \, \sum_{|x|=k} {\bf 1}_{\{ x< \mathscr{L}_n\} } \, \ee^{-V(x)} 
\le
\sum_{k=b(\log n)^2}^{B(\log n)^2} W_k^{(\log n)} ,
$$ 

\noindent and 
$$
\sum_{k=b(\log n)^2}^{B(\log n)^2} \, \sum_{|x|=k} {\bf 1}_{\{ x< \mathscr{L}_n\} } \, \ee^{-V(x)} 
\ge
\sum_{k=b(\log n)^2}^{B(\log n)^2} W_k^{(\log\frac{n}{B (\log n)^2})} .
$$

\noindent Applying \eqref{Zn-proof5} to $\lambda = \log n$ and noting that $\lim_{n\to \infty} \frac{\log\frac{n}{B (\log n)^2}}{\log n} = 1$, this yields that for any fixed $B>b>0$,
\begin{equation}
    I_{\eqref{pnbB}}  
    =
    \big(\frac{8\sigma^2}{\pi}\big)^{1/2} \, \e \Big[ [(B^{1/2} \wedge \frac{1}{\sigma \eta}) - b^{1/2}]\, {\bf 1}_{( \eta \le \frac{1}{\sigma \, b^{1/2}})}\Big] 
    +
    o_{\P^*}(1) \, .
    \label{pnbB2}
\end{equation}
  
\noindent Note that $\e \{ [(B^{1/2} \wedge \frac{1}{\sigma \eta}) - b^{1/2}]\, {\bf 1}_{( \eta \le \frac{1}{\sigma \, b^{1/2}})}\}$ is continuous in $B$ and $b$. Since $|X_n|$ and $|X_{n-1}|$ only differ $1$, we get that 
\begin{eqnarray*}
 &&2 P_\omega  \Big( b + \frac{1}{(\log n)^2} \le \frac{|X_n|}{(\log n)^2} \le B - \frac{1}{(\log n)^2} \Big)  
    \\
 &\le& I_{\eqref{pnbB}}  
    \le 
    2 P_\omega  \Big( b- \frac{1}{(\log n)^2} \le \frac{|X_n|}{(\log n)^2} \le B + \frac{1}{(\log n)^2} \Big) ,
\end{eqnarray*}

\noindent which, in view of \eqref{pnbB2}, readily yields the case $\kappa=1$ of Corollary \ref{c:cv}, as claimed.\hfill$\Box$

\bigskip

\noindent {\bf Acknowledgements.} We are grateful to Dayue Chen, Nina Gantert and Ofer Zeitouni for stimulating discussions on biased walks, and to Laurent Saloff-Coste for enlightenment on Markov chains.

\end{document}